\documentclass[11pt]{amsart}\usepackage{CJK}
\usepackage[centertags]{amsmath}
\usepackage{times,fullpage}
\usepackage{latexsym,amscd,amssymb}
\usepackage{tikz,tikz-cd}
\usepackage{dutchcal}
\usepackage[all]{xy}
 \input xy
  \xyoption{all}

\usepackage{amsfonts}
\usepackage{amssymb}
\usepackage{amsthm}
\usepackage{newlfont}
\usepackage[colorlinks,linkcolor=blue, pdfstartview=FitH]{hyperref}
\usepackage{mathrsfs}
\usepackage{enumerate}
\numberwithin{equation}{section}

 \theoremstyle{plain}
\newtheorem{thm}{Theorem}[section]
\newtheorem{theorem}[thm]{Theorem}
\newtheorem{lemma}[thm]{Lemma}
\newtheorem{corollary}[thm]{Corollary}
\newtheorem{proposition}[thm]{Proposition}

\theoremstyle{definition}

\newtheorem{setting}[thm]{Setting}

\newtheorem{question}[thm]{Question}
\newtheorem{conjecture}[thm]{Conjecture}

\newtheorem{remark}[thm]{Remark}

\newtheorem{definition}[thm]{Definition}
\newtheorem{claim}[thm]{Claim}

\newtheorem{example}[thm]{Example}

\newtheorem{defn-thm}[thm]{Definition-Theorem}

\newcommand{\sE}{{\mathcal E}}
\newcommand{\sF}{{\mathcal F}}
\newcommand{\sG}{{\mathcal G}}

\newcommand{\sI}{{\mathcal I}}

\newcommand{\sL}{{\mathcal L}}

\newcommand{\sO}{{\mathcal O}}

\newcommand{\sQ}{{\mathcal Q}}

\newcommand{\A}{{\mathbb A}}

\newcommand{\C}{{\mathbb C}}

\newcommand{\N}{{\mathbb N}}
\renewcommand{\P}{{\mathbb P}}
\newcommand{\Q}{{\mathbb Q}}
\newcommand{\R}{{\mathbb R}}

\newcommand{\Z}{{\mathbb Z}}

\newcommand{\btheorem}{\begin{theorem}}
\newcommand{\etheorem}{\end{theorem}}
\newcommand{\bquestion}{\begin{question}}
\newcommand{\equestion}{\end{question}}
\newcommand{\bconjecture}{\begin{conjecture}}
\newcommand{\econjecture}{\end{conjecture}}
\newcommand{\bclaim}{\begin{claim}}
\newcommand{\eclaim}{\end{claim}}
\newcommand{\bsetting}{\begin{setting}}
\newcommand{\esetting}{\end{setting}}
\newcommand{\bproposition}{\begin{proposition}}
\newcommand{\eproposition}{\end{proposition}}
\newcommand{\bdefinition}{\begin{definition}}
\newcommand{\edefinition}{\end{definition}}
\newcommand{\bcorollary}{\begin{corollary}}
\newcommand{\ecorollary}{\end{corollary}}
\newcommand{\bproof}{\begin{proof}}
\newcommand{\eproof}{\end{proof}}
\newcommand{\bremark}{\begin{remark}}
\newcommand{\eremark}{\end{remark}}
\newcommand{\eexample}{\end{example}}
\newcommand{\bexample}{\begin{example}}
\newcommand{\elemma}{\end{lemma}}
\newcommand{\blemma}{\begin{lemma}}

\newcommand{\emb}{\hookrightarrow}

\renewcommand{\bar}{\overline}

\renewcommand{\phi}{\varphi}

\newcommand{\ee}{\end{eqnarray*}}
\newcommand{\be}{\begin{eqnarray*}}

\newcommand{\beq}{\begin{equation}}
\newcommand{\eeq}{\end{equation}}

\newcommand{\bd}{\begin{enumerate}[(1)]}
\newcommand{\ed}{\end{enumerate}}

\renewcommand{\hat}{\widehat}
\renewcommand{\tilde}{\widetilde}

\renewcommand{\bf}{\textbf}

\renewcommand{\rm}{\textrm}

\address{Beijing International Center for Mathematical Research,
Peking University, No. 5 Yiheyuan Road, Haidian District, Beijing 100871, China}
\email{zhengxu@pku.edu.cn}


\begin{document}

\title{Abundance for threefolds in positive characteristic when $\nu=2$}
\author{Zheng Xu}

\date{}

\maketitle

\begin{abstract}
In this paper, we prove the abundance conjecture for threefolds over a perfect field $k$ of characteristic $p > 3$ in the case where the numerical dimension is equal to $2$. More  precisely, we show that if $(X,B)$ is a projective log canonical pair of dimension $3$ over $k$ such that $K_{X}+B$ is nef and $\nu(K_{X}+B)=2$, then $K_{X}+B$ is semi-ample. 

\end{abstract}

{\tableofcontents \setcounter{page}{0}\pagenumbering{roman}}

\setcounter{page}{1}\pagenumbering{arabic}

\section{Introduction}

In birational geometry, the Minimal Model Program (MMP) is a conjectural program to classify all algebraic varieties up to birational equivalence. This program predicts that after finitely many birational transforms, every variety with mild singularities admits a birational model which is either a minimal model or a Mori fibre space. In characteristic $0$, much progress has been made towards the MMP. For example, it has been fully established in dimensions $\leq 3$. In higher dimensions, it's proved in \cite{birkar2010existence} that the minimal models for varieties of general type exist.
However, in positive characteristics, due to the failure of vanishing theorems, there has been some big progress in the MMP only after the work of Hacon and Xu (see \cite{hacon2015three}).
They proved the existence of minimal models for terminal threefolds over an algebraically closed field $k$ of characteristic $p > 5$. Then Cascini, Tanaka and Xu proved that arbitrary terminal threefold over $k$ is birational to either a minimal model or a Mori fibre space (\cite{cascini2015base}). Based on it, Birkar and Waldron established
the MMP for klt threefolds over $k$ (\cite{birkar2016existence,birkar2017existence}). Moreover, there are some generalizations of it in various directions. For example, see \cite{hashizume2020minimal,waldron2018lmmp} for its generalization to log canonical (lc for short) pairs, \cite{gongyo2019rational,hacon2022minimal,hacon2022relative} for its generalization to low characteristics, \cite{das2022log,waldron2023mori} for its generalization to imperfect base fields, and \cite{bhatt2023globally,takamatsu2023m} for its analog in mixed characteristics.

Now we can run MMP for lc pairs of dimension $3$ over a perfect field of characteristic $p >3$ (see Theorem \ref{lcmmp}). Hence a central problem remaining is the following conjecture.\\

\noindent \bf{Abundance conjecture.} 
Let $(X,B)$ be a projective lc pair of dimension $3$ over a perfect field $k$ of characteristic $p >3$. If $K_{X}+B$ is nef,
then it is semi-ample.\\

In characteristic $0$, the abundance conjecture for threefolds was  proved in \cite{kawamata1992abundance,keel1994log}. The proof has to be parted according to the numerical dimension $\nu(K_{X}+B)$. In characteristic $p>3$, the conjecture has been proved in the cases of  $\nu(K_{X}+B)=0,3$. When $\nu(K_{X}+B)=1,2$, it is open even in the terminal case. 
For details, please see \cite{birkar2017existence,das2019abundance,hacon2022minimal,waldron2017finite,waldron2018lmmp,
witaszek2021canonical,xu2024note,zhang2019abundance,zhang2020abundance}
for example. In mixed characteristic, the abundance conjecture for arithmetic threefolds whose closed points have residue characteristic $>5$ was proved in \cite{bernasconi2024abundance}.
In this paper, we prove the conjecture in the case of $\nu(K_{X}+B)=2$.

\btheorem (Theorem \ref{mainthm}) Let $(X,B)$ be a projective lc pair of dimension $3$ over a perfect field $k$ of characteristic $p>3$. If $K_{X}+B$ is nef and $\nu(K_{X}+B)=2$, then it is semi-ample.  
\etheorem

When $k$ is $\C$, this result is proved in the terminal case in \cite{kawamata1992abundance} (see \cite[Theorem 14.4.1]{kollar1992flips} for a simplified proof). It is proved in the lc case in \cite{keel1994log}. The proofs rely on numerous results on the intersection theory for $\Q$-varieties together with positivity theorems, such as vanishing theorems and generic semi-positivity theorems. These results may fail in positive characteristics. \\

\noindent \bf{Sketch of the proof of Theorem 1.1.}

\noindent For simplicity, we assume that $k$ is an algebraically closed field of characteristic $p>5$ since in this case, klt singularities of dimension $3$ are rational singularities by \cite[Corollary 1.3]{arvidsson2022kawamata}.  By a standard lemma (see Lemma \ref{standardmodification}), we can reduce the assertion to the case when $(X,B)$ is an lc projective pair of dimension $3$ and the following assertions hold, 

\noindent (1) $X$ is $\Q$-factorial, $B$ is reduced, and $X\backslash B$ is terminal,

\noindent (2) $K_{X}+B\sim_{\Q} D$ for an effective $\Q$-divisor $D$ with $\mathrm{Supp}\ D=  B,$

\noindent (3) $K_{X}+(1-\varepsilon)B$ is nef for any sufficiently small $\varepsilon>0$,

\noindent (4)  $\nu(K_{X}+B)=2$,

\noindent (5) if $C$ is a curve in $X$ with $C\cdot (K_{X}+B)>0$, then $(X,B)$ is plt at the generic  point of $C$.

\noindent In this special case, we can mimic the proof in characteristic $0$ (see \cite[Theorem 14.4.1]{kollar1992flips}). Let $\rho:V\to X$ be a desingularization of $X$ and $m$ be a positive integer such that $L:=m(K_{X}+B)$ is Cartier. By Riemann-Roch theorem, we can prove that
$$\chi(X,\sO_{X}(nL))=\frac{n}{12}(K^{2}_{V}+c_{2}(V))\cdot \rho^{\ast}L+\chi(\sO_{V}).$$
Here we use the assumption of the characteristic so that $X$ has rational singularities.
Then we meet two cases.

\noindent Case \uppercase\expandafter{\romannumeral1}: $T_{V}$ is strongly $(\rho^{\ast}L,\rho^{\ast}L)$-semistable (see Definition \ref{stabilitydefinition}).

\noindent Case \uppercase\expandafter{\romannumeral2}: $T_{V}$ is not strongly $(\rho^{\ast}L,\rho^{\ast}L)$-semistable.

In Case \uppercase\expandafter{\romannumeral2}, we use Langer's results on bend and break by $1$-foliations (see Theorem \ref{bendbreakfoliation}) to prove that the nef dimension $n(K_{X}+B)\leq 2$. Then the assertion follows from Theorem \ref{known}.
Case \uppercase\expandafter{\romannumeral1} is much more difficult since we don't have a good intersection theory for singular varieties in positive characteristic. Instead of the intersection theory for $\Q$-sheaves used in characteristic $0$, we use the intersection theory developed by Langer in \cite{langer2025intersection}.
In this case, $\Omega_{X}^{[1]}:=(\Omega_{X})^{\ast\ast}$ is also strongly $(L,L)$-semistable since $\rho$ is a birational morphism.
Then we can apply Langer's Bogomolov inequality for reflexive sheaves (see Theorem \ref{Bogomolovineq}) to prove that 
$$L\cdot\big(K_{X}^{2}+c_{2}(\Omega_{X}^{[1]})\big)\geq 0,$$
where the second Chern class is defined in Section $4$.
Then the most important step of the proof is to show the following key inequality
$$\rho^{\ast}L \cdot\big(K^{2}_{V}+c_{2}(V)\big)\geq L\cdot\big(K_{X}^{2}+c_{2}(\Omega_{X}^{[1]})\big).$$
To do this, we prove some Bertini type results in characteristic $p>3$ (see Section 3) and
reduce the proof of the key inequality to a problem on surfaces (see Proposition \ref{keyinequality}).
Then we need to compute relative Chern classes (see Section 4 for definition) of desingularizations of some cyclic quotient surfaces singularities. We explain how to compute these classes in Subsection \ref{subsection 5.4}.
Therefore, we can prove the key inequality and modify the proof of \cite[Theorem 14.4.1]{kollar1992flips} to prove that $\kappa(K_{X}+B)> 0$, and hence the assertion follows from Theorem \ref{known}. \\

\noindent\bf{Notation and conventions.}

$\bullet$ All schemes considered in this article are excellent, in particular Noetherian.

$\bullet$ A singular point of a scheme is a point that is not regular.

$\bullet$ A variety is an integral and separated scheme which is  of finite type over a field $k$.

$\bullet$ We call $(X,B)$ a pair if $X$ is a normal variety and $B$ is an effective $\Q$-divisor on $X$ such that $K_{X}+B$ is $\Q$-Cartier. For more notions in the theory of MMP 
such as klt (dlt, lc, plt) pairs, filps, divisorial contractions and so on, we refer to \cite{kollar1998birational}.

$\bullet$ Let $X$ be a normal projective variety over a field $k$ and $D$ be a $\Q$-Cartier $\Q$-divisor on $X$. If $|mD|=\emptyset$ for all $m>0$, we define the Kodaira dimension $\kappa(X,D)=-\infty$. Otherwise, let $\Phi:X\dashrightarrow Z$ be the Iitaka map (we refer to \cite[2.1.C]{lazarsfeld2017positivity}) of $D$ and we define the Kodaira dimension $\kappa(X,D)$ to be the dimension of the image of $\Phi$. Sometimes, we write $\kappa(D)$ for $\kappa(X,D)$. We denote $\kappa(X,K_{X})$ by $\kappa(X)$. For a projective variety $Y$ over a field $k$ admitting a smooth model $\tilde{Y}$, we define $\kappa(Y):=\kappa(\tilde{Y})$.

$\bullet$ Let $X$ be a normal projective variety of dimension $n$ over a field $k$ and $D$ be a   nef $\Q$-Cartier $\Q$-divisor on $X$.
Then we can define
$$\nu(X,D):=\mathrm{max}\{k\in\N | D^{k}\cdot A^{n-k}>0\ \rm{for}\ \rm{an}\ \rm{ample}\ \rm{divisor}\ A\ \rm{on}\ X\}.$$
Sometimes, we write $\nu(D)$ for $\nu(X,D)$.

$\bullet$ Assume that $f:X\to Y$ is a morphism between normal schemes, $\sF$ is a reflexive sheaf on $X$ and $\sE$ is a reflexive sheaf on $Y$. Then we set
$$f^{[\ast]}\sE:=(f^{\ast}\sE)^{\ast\ast},f_{[\ast]}\sF:=(f_{\ast}\sF)^{\ast\ast}.$$ 
\\

\noindent\bf{Acknowledgements.}
I would like to express my gratitude to my advisor Wenhao Ou for his help,
encouragement, and support. I would like to thank Lei Zhang for his encouragement and advice. I thank the anonymous reviewers for many corrections and suggestions. The author is supported by the National Key R$\&$D Program of China (No.2021YFA1002300).

\section{Preliminaries}

In this section we recall some basic results.

\subsection{Nef reduction map}

In this subsection, we recall the notion of nef reduction map.

\bdefinition Let $X$ be a normal projective variety defined over an uncountable field
 and let $L$ be a nef $\Q$-Cartier  $\Q$-divisor. We call a rational map 
$\phi: X \dashrightarrow Z$ a \textit{nef reduction map} of $L$ if $Z$ is a normal projective variety and there exist open dense subsets $U \subseteq X$, $V\subseteq Z$ such that

\noindent (1)  $\phi|_{U}: U \to V$ is proper  and 
$\phi_{\ast}\sO_{U} = \sO_{V}$ ,

\noindent (2)  $L|_{F}\equiv 0$ for all fibres $F$ of $\phi$ over $V$, and

\noindent (3)  if $x \in X$ is a very general point and $C$ is a curve passing through it, then $C\cdot L = 0$ if and only if $C$ is contracted by $\phi$. 
\edefinition

It is proved that the nef reduction map exists over an uncountable algebraically closed field.

\btheorem (\cite[Theorem 2.1]{bauer2002reduction}) A nef reduction map exists for normal projective varieties defined over an uncountable algebraically closed field. Further, it is unique
up to birational equivalence.
\etheorem

We call $n(X,L)= \mathrm{dim}\ Z$ the \textit{nef dimension} of $L$, where $Z$ is the target of a nef reduction map $\phi: X \dashrightarrow Z$. When the base field is countable and algebraically closed, we can define $n(X,L):=n(X_{K},L_{K})$ by \cite[Proposition 2.16]{witaszek2021canonical}, where $K$ is an uncountable algebraically closed field that contains $k$, and $X_{K},L_{K}$ are the base changes of $X,L$ to $K$.
It satisfies that $\kappa(X, L)\leq n(X, L)$.

\subsection{Surface singularities and slc surfaces}

In this subsection, we recall the notion of the dual graph of a surface singularity and collect some results on lc surface singularities and slc surfaces. 
The main references for this subsection are \cite{kollar2013singularities,sato2025general}.
First we recall some basic knowledge on the intersection theory for curves over a  field, which is not necessarily perfect. Let $C$ be an integral projective scheme over a field $k$ with $\mathrm{dim}\ C =1$ (not necessarily regular).

\bdefinition\label{genusdefinition} We define the \textit{arithmetic genus} $g(C)$ of $C$ by 
$$g(C):=\frac{\mathrm{dim}_{k}(H^{1}(C,\sO_{C}))}{\mathrm{dim}_{k}(H^{0}(C,\sO_{C}))}.$$
For a $0$-cycle $D =\sum_{i=1}^{n} a_{i}P_{i}$ on $C$, we define the \textit{degree} of $D$ over $k$ by
$$\mathrm{deg}_{C/k}(D):=\sum_{i=1}^{n}a_{i}\ \mathrm{dim}_{k}(k(P_{i}))\in \Z$$
where $k(P_{i})$ are the residue fields of $C$ at $P_{i}$. Note that any line bundle $L\in \mathrm{Pic}(C)$ is isomorphic to $\sO_C(B)$ for some $0$-cycle $B$ on $C$. We define 
$$\mathrm{deg}_{C/k}(L):=\mathrm{deg}_{C/k}(B).$$ 
Since rationally equivalent cycles have the same degree, $\mathrm{deg}_{C/k}(L)$ is well-defined.
The map $\mathrm{deg}_{C/k}$ gives the homomorphism 
$$\mathrm{deg}_{C/k}:\mathrm{Pic}(C)\to \Z.$$
See \cite[Subsection 1.4]{fulton2013intersection} for details. 
\edefinition

Let $X$ be a normal integral scheme with a dualizing complex $\omega_{X}^{\bullet}$. Then we can define resolutions of singularities (desingularizations for short), the canonical divisor $K_{X}$, and lc (klt, plt) pairs as usual (see \cite[Section 2]{sato2025general} for example).

\bdefinition We say that $(X,x)$ (or $(x\in X)$) is a \textit{normal surface singularity} if $X$
is a two-dimensional normal integral scheme with a dualizing complex $\omega_{X}^{\bullet}$ and $x$ is a closed point of $X$.
\edefinition
\bremark The definition is more general than that in \cite[Definition 2.10]{sato2025general}, as it does not require that $X$ be the spectrum of a local ring.
For an effective $\Q$-Weil divisor $B$ on $X$, we also say that $(x\in X,B)$ is a normal surface singularity.
\eremark

Let $(X,x)$ be a normal surface singularity, and let $f: Y\to X$ be a proper birational morphism from a normal integral scheme $Y$ with the exceptional divisor $\mathrm{Exc}(f)= \sum_{i}E_{i}$ satisfying $f(\mathrm{Exc}(f))=x$. Note that $E_{i}$ are defined over the residue field $k(x)$ of $X$ at $x$.

\bdefinition\label{intersectionnumber} For a Cartier divisor $D$ on $Y$ and an exceptional Weil divisor $Z =\sum_{i}b_{i}E_{i}$ on $Y$, we define the \textit{intersection number} $D \cdot Z$ by
$$D \cdot Z:=\sum_{i}b_{i}\ \mathrm{deg}_{E_{i}/k(x)}(\sO_{Y}(D)|_{E_{i}}).$$
\edefinition

By \cite[Definition 2.11 and Remark 2.12]{sato2025general}, the \textit{minimal resolution} of a normal surface singularity $(X,x)$ can be defined, and it exists uniquely.
Now we can give the definition of the dual graphs as follows.

\bdefinition\label{dualgraphdefinition}
Let $(X,x)$ be a normal surface singularity, and $f: Y\to X$ be the minimal resolution with the exceptional divisor $\mathrm{Exc}(f) = \sum_{i=1}^{n} E_{i}$ satisfying $f(\mathrm{Exc}(f))=x$.
The \textit{dual graph} of $(X,x)$ is a graph whose set of vertices is $\{E_{1}, \cdots , E_{n}\}$, the
number of edges between $E_{i}$ and $E_{j}$ is $(E_{i}\cdot E_{j} )\in \N$ and the weight at $E_{i}$ is
$$(\mathrm{dim}_{k(x)} H^{0}(E_{i}, \sO_{E_{i}}), g(E_{i}),E_{i}^{2}
)\in \Z^{3}.$$
\edefinition
\bremark We can define the \textit{extended dual graph} of a normal surface singularity $(x\in X,B)$ similarly (see \cite[2.26]{kollar2013singularities} for details). We often omit the weights when they are either unknown or not relevant. When we work with klt surface singularities over algebraically closed fields, the notion of dual graphs can be largely simplified since $\mathrm{dim}_{k(x)} H^{0}(E_{i}, \sO_{E_{i}})=1$ and $g(E_{i})=0$ (see \cite[Theorem 4.7]{kollar1998birational} for example).
\eremark

\bdefinition (\cite[Corollary 3.31]{kollar2013singularities})\label{definition_cyclic quotient singularities} Let $S$ be a two-dimensional normal integral scheme with a dualizing complex, $s\in S$ a closed point and $B$ a reduced Weil divisor on $S$. Assume that $(s\in S,B)$ is lc and let $(\Gamma,B)$ be its extended dual graph. We denote the residue field of $S$ at $s$ by $k(s)$. We say $(s\in S,B)$ is a \textit{cyclic quotient singularity} if $(\Gamma,B)$ is
$$\bullet\rule[3pt]{1cm}{0.05em}\circ\rule[3pt]{1cm}{0.05em}\cdots\rule[3pt]{1cm}{0.05em}\circ,\ \mathrm{or}\ \circ\rule[3pt]{1cm}{0.05em}\circ\rule[3pt]{1cm}{0.05em}\cdots\rule[3pt]{1cm}{0.05em}\circ,$$ 
where the solid bullet stands for the strict transform of $B$ and every circle stands for an exceptional curve over $S$. Moreover, every exceptional curve is a $\P^1_{k(s)}$, except in the case $B=0$ with only one exceptional curve, which may instead be a $k(s)$-irreducible conic.
\edefinition

We will use the following result in Corollary \ref{bertiniforcyclicquotient} and Section 5.

\bproposition\label{lcsurfacesing}  Let $S$ be a two-dimensional normal integral scheme with a dualizing complex, $s\in S$ a closed point and $B$ a reduced Weil divisor on $S$. Assume that $(s\in S,B)$ is lc and let $(\Gamma,B)$ be its extended dual graph. If $(s\in S,B)$ is plt and $s\in B\neq 0$,
then $(\Gamma,B)$ is 
$$\bullet\rule[3pt]{1cm}{0.05em}\circ\rule[3pt]{1cm}{0.05em}\cdots\rule[3pt]{1cm}{0.05em}\circ,$$ 
where the solid bullet stands for the strict transform of $B$, every circle stands for $\P^{1}_{k(s)}$ and $k(s)$ is the residue field of $S$ at $s$. In particular, $(s\in S,B)$ is a cyclic quotient singularity.
\bproof By \cite[Corollary 3.31]{kollar2013singularities}, the pair $(s\in S,B)$ must be one of the follow singularities: (1), (2) or (3) in \cite[Corollary 3.31]{kollar2013singularities}, or those in \cite[Examples 3.27 and 3.28]{kollar2013singularities}. Singularities in (2) and (3) of \cite[Corollary 3.31]{kollar2013singularities}, as well as those in \cite[Example 3.27]{kollar2013singularities}, require that $B=0$. Singularities in \cite[Example 3.28]{kollar2013singularities} require that $(s\in S,B)$ is not plt. Hence we are in case (1) of \cite[Corollary 3.31]{kollar2013singularities}. By Definition \ref{definition_cyclic quotient singularities}, this means that $(s\in S,B)$ is a cyclic quotient singularity. Moreover, since $B\neq 0$, the assertion holds.
\eproof\eproposition

We recall the following adjunction formula of lc surface singularities.

\btheorem\label{lcadjunction} (\cite[Theorem 3.36]{kollar2013singularities}) Let $S$ be a two-dimensional normal integral scheme with a dualizing complex, and let $B\subset S$ be a one-dimensional reduced proper scheme over a field  such that $(S,B)$ is lc. Then
$$(K_{S}+B)\cdot B=\mathrm{deg}\ \omega_{B}+\sum_{s\ \mathrm{cyclic,\ plt}}(1-\frac{1}{\mathrm{det}(\Gamma_{s})})\mathrm{deg}[s]+\sum_{s\ \mathrm{dihedral}}\mathrm{deg}[s],$$
where $\mathrm{det}(\Gamma_{s})$ are the determinants of the dual graphs $\Gamma_{s}$ of cyclic quotient singularities $(S,s)$.
\etheorem
\bremark\label{cyclic_quotient_action} We also call $\mathrm{det}(\Gamma_{s})$ the index of $(S,s)$. 
By \cite[Theorem 3.32]{kollar2013singularities} for a cyclic quotient singularity $(S,s)$,
if the residue field $k(s)$ of $S$ at $s$ is algebraically closed and $S$ is defined over $k(s)$, then the completion of $S$ at $s$ is isomorphic to 
$$\hat{\mathbb{A}}_{k(s)}^{2}/\mu_{m},$$
where $\mu_{m}$ is the group scheme $\mathrm{Spec}(k(s)[\Z_{m}])$. By \cite[Definition 3.33]{kollar2013singularities}, we have $\mathrm{det}(\Gamma_{s})=m$.  If moreover $m$ is not divisible by the characteristic of $k(s)$, then by \cite[3.19]{kollar2013singularities}, the action of $\mu_{m}\cong \Z_{m}$ on $\hat{\mathbb{A}}_{k(s)}^{2}$ is linear and diagonalized.
\eremark

Now we turn to the notion of slc surfaces.

\bdefinition Let $X$ be a scheme. A coherent $\sO_{X}$-module $\sF$ is $S_{i}$ where $i$ is a positive integer if it satisfies the condition
$$\mathrm{depth}_{\sO_{X,x} \sF_{x}} \geq \mathrm{min}\{i,\mathrm{dim}(\mathrm{Supp}\ \sF_{x})\}$$ for all $x \in X$. We say that $X$ is $S_{i}$ if $\sO_{X}$ is $S_{i}$.
\edefinition

\bproposition\label{S2fication} (\cite[Proposition 2.6]{waldron2018lmmp}) Let $X$ be a reduced equidimensional quasi-projective scheme. Then
the locus $U$ where $X$ is $S_{2}$ is an open subset with $\mathrm{codim}_{X} (X\backslash U)\geq 2$, and there exists a birational morphism called $S_{2}$-fication $g: X^{\prime}\to X$ satisfying:

\noindent (1) $X^{\prime}$ is $S_{2}$ and reduced,

\noindent (2) $g$ is finite and an isomorphism precisely above $U$, and

\noindent (3) the normalization $X^{n}\to X$ factorizes through $g$.
\eproposition

\bdefinition A scheme (or ring) is called \textit{nodal} if its codimension one local
rings are regular or nodal. It is called \textit{demi-normal} if it is $S_{2}$ and nodal.
\edefinition

We recall the following Riemann-Roch theorem for demi-normal surfaces.

\btheorem\label{RRslc} Let $X$ be a demi-normal surface over an algebraically closed field. Let $L$ be a Cartier divisor on $X$. Then we have
$$\chi(\sO_{X}(L))=\frac{1}{2}L\cdot(L-K_{X})+\chi(\sO_{X}).$$ 
\bproof In characteristic $0$, it is proved in \cite[Theorem 3.1]{liu2016geography}. 
Assume that $X$ is defined over an algebraically closed field of arbitrary characteristic. Let $\pi:\overline{X}\to X$ be the normalization of $X$, and let $D\subseteq X$ and $\overline D\subseteq \overline{X}$ be conductor subschemes.  By the proof of \cite[Theorem 3.1]{liu2016geography}, it suffices to show that 
$$\chi(\sO_X)=\chi(\sO_{\overline X})+\chi(\sO_D)-\chi(\sO_{\overline{D}})$$
and 
$$\chi(\sO_X(L))=\chi(\pi^\ast\sO_X(L))+\chi(\sO_X(L)_{|_D})-\chi(\pi^\ast\sO_X(L)_{|_{\overline{D}}}).$$
We consider the exact sequences
$$0\to \sO_{\overline X}(-\overline{D})\to \sO_{\overline X}\to \sO_{\overline{D}}\to 0$$
and
$$0\to \sI_{D}\to \sO_X\to \sO_{D}\to 0.$$
By the definition of conductor subschemes \cite[5.2]{kollar2013singularities}, we have
$\pi_\ast \sO_{\overline X}(-\overline{D})\cong \sI_{D}.$ It follows that
$$\chi(\sO_X)-\chi(\sO_D)=\chi(\sI_{D})=\chi(\pi_\ast \sO_{\overline X}(-\overline{D}))=\chi(\sO_{\overline X})-\chi(\sO_{\overline{D}}).$$
Similarly, we have
$$\chi(\sO_X(L))-\chi(\sO_X(L)_{|_D})=\chi(\pi^\ast\sO_X(L))-\chi(\pi^\ast\sO_X(L)_{|_{\overline{D}}}).$$
Therefore, the assertion holds.
\eproof
\etheorem

\bdefinition (\cite[5.10]{kollar2013singularities}) We say that $(X,\Delta)$ is a \textit{semi-log canonical (slc)} pair if: $X$ is demi-normal, $\Delta$ is an effective $\Q$-divisor with no components along $D$, $K_{X}+\Delta$ is $\Q$-Cartier, and the normalization $(\overline{X},\overline{D}+\overline{\Delta})$ is an lc pair, where $D$ and $\overline{D}$ are conductor subschemes, and $\overline{\Delta}$ is the birational transform of $\Delta$. \edefinition

The abundance for slc surfaces over any field of positive characteristic holds. 

\btheorem\label{slcabundance} (\cite[Theorem 1]{posva2023abundance}) Let $(X,\Delta)$ be a projective slc surface over a field of characteristic $p >0$. If $K_{X}+\Delta$ is nef, then $K_{X}+\Delta$ is semi-ample.
\etheorem

In characteristic $0$, the reduced part of the boundary of an lc pair of dimension $3$ is slc. It fails in positive characteristic, despite several partial positive results (see \cite{arvidsson2023vanishing,arvidsson2023normality,arvidsson2024properness} for example). However, we have the following analogous result.

\blemma\label{slcsurface} Let $(X,B)$ be an lc projective pair of dimension $3$ over an algebraically closed field $k$ of characteristic $p>0$. Assume that $B$ is reduced and let $g:B^{\prime}\to B$ be its $S_{2}$-fication.
Then the pair $(B^{\prime},\mathrm{Diff}_{B^{\prime}}(0))$ defined by the adjunction $K_{B^{\prime}}+\mathrm{Diff}_{B^{\prime}}(0)=(K_{X}+B)|_{B^{\prime}}$ is a slc pair, where the existence of $\mathrm{Diff}_{B^{\prime}}(0)$ is given by \cite[Definition 4.2]{kollar2013singularities}.
\bproof By the classification of lc surface singularities (see \cite[4.1]{kollar1998birational} for example), we know $B$ is nodal. Since $g$ is an isomorphism outside some points by Proposition \ref{S2fication}, we have $B^{\prime}$ is nodal. Therefore, $B^{\prime}$ is demi-normal. Then by \cite[Theorem 17.2]{kollar1992flips} we know $(B^{\prime},\mathrm{Diff}_{B^{\prime}}(0))$ is slc.
\eproof
 \elemma

\subsection{MMP for threefolds in positive characteristic}

In this subsection, we recall the theory of MMP for projective lc pairs of dimension $3$ over a perfect field of characteristic $p> 3$. Moreover, we define a partial MMP over algebraically closed fields of characteristic $p> 3$ (see Definition \ref{MMPKtrivial}) and we will use this construction to study the abundance in Section $5$.

\btheorem\label{lcmmp} (\cite[Theorem 1.1]{hashizume2020minimal} and \cite{hacon2022minimal})  Let $(X, B)$ be an lc pair of dimension $3$ over a perfect field $k$ of characteristic $p>3$, and let $f : X \to Y$ be a projective surjective morphism to a quasi-projective variety. If $K_{X}+B$ is pseudo-effective (resp. not pseudo-effective) over $Y$ , then we can run a $(K_{X}+B)$-MMP to get a log minimal model (resp. Mori fibre space) over $Y$ .
\etheorem

We recall the notion of MMP with scaling.
Let $(X,B)$ be a projective lc pair of dimension $3$ over a perfect field $k$ of characteristic $p>3$ and $A\neq 0$ an effective $\Q$-Cartier $\Q$-divisor on $X$. Suppose that there is $t_{0}> 0$ such that $(X,B+ t_{0}A)$ is lc and $K_{X}+B+t_{0}A$ is nef. We describe how to run a $(K_{X}+B)$-MMP with scaling of $A$ as follows.

Let $\lambda_{0}=\mathrm{inf} \{t: K_{X}+B+tA\ \mathrm{is}\ \mathrm{nef} \}$. Suppose we can find a $(K_{X}+ B)$-negative extremal ray $R_{0}$ which satisfies $(K_{X}+B+\lambda_{0}A)\cdot R_{0}= 0$ (In general, it is possible that there is no such extremal ray). This is the first ray we contract in our MMP. If the contraction
is a Mori fibre contraction, we stop. Otherwise let $X_{1}$ be the result of the
divisorial contraction or flip. Then $K_{X_{1}}+ B_{X_{1}}+\lambda_{0}A_{X_{1}}$
is also nef, where $B_{X_{1}}$ and $A_{X_{1}}$ denote the birational transforms on $X_{1}$ of $B$ and $A$, respectively. We define $\lambda_{1} =\mathrm{inf} \{t: K_{X_{1}}+ B_{X_{1}} + tA_{X_{1}}\
\mathrm{is}\ \mathrm{nef} \}$. The next step in our MMP is chosen to be a $(K_{X_{1}} + B_{X_{1}})$-negative extremal ray $R_{1}$ which is $(K_{X_{1}}+B_{X_{1}}+ \lambda_{1} A_{X_{1}})$
-trivial. So long as we can always find the appropriate extremal
rays, contractions and flips, we can continue this process.

\bproposition\label{findextremalray} Let $(X,B)$ be a $\Q$-factorial projective lc pair of dimension $3$ over an algebraically closed field $k$ of characteristic $p>3$, and let $W$ be an effective $\Q$-divisor such that
$K_{X}+B+W$ is nef. Then either

\noindent (1) there is a $(K_{X}+B)$-negative extremal ray which is $(K_{X}+B+W)$-trivial, or

\noindent (2) $K_{X}+B+(1-\varepsilon)W$ is nef for any sufficiently small rational $\varepsilon>0$.
\bproof It is an adaptation of \cite[Lemma 5.1]{keel1994log}. Note that the proof there only uses the fact that for any $(K_{X}+B)$-negative extremal ray $R$ there is a rational curve $C$ such that $C$ generates $R$ and $-(K_{X}+B)\cdot C\leq 6$, which holds in our setting by \cite[Theorem 1.3]{hashizume2020minimal} and \cite{hacon2022minimal} since $k$ is an algebraically closed field of characteristic $p>3$.
\eproof
\eproposition

\bcorollary\label{MMPwithscaling} Let $(X,B)$ be a $\Q$-factorial projective lc pair of dimension $3$ over an algebraically closed field $k$ of characteristic $p>3$ and $A$ be an effective $\Q$-divisor such that $(X,B+A)$ is lc and $K_{X}+B+A$ is nef.
If $K_{X}+B$ is not nef, then we can run a $(K_{X}+B)$-MMP with scaling of $A$.
\bproof Let $\lambda:= \mathrm{inf} \{t: K_{X}+B+tA\ \mathrm{is}\ \mathrm{nef} \}$. It suffices to show that we can find a $(K_{X}+B)$-negative extremal ray $R$ such that $(K_{X}+B+\lambda A)\cdot R=0$. We apply Proposition \ref{findextremalray} by letting $W:=\lambda A$.
\eproof
\ecorollary

In this paper, we will use the following construction.

\bdefinition\label{MMPKtrivial} Let $(X,B)$ be a $\Q$-factorial projective lc pair of dimension $3$ over an algebraically closed field $k$ of characteristic $p>3$, and let $A$ be an effective $\Q$-divisor such that $(X,B+A)$ is lc and $K_{X}+B+A$ is nef. We can run a partial $(K_{X}+B)$-MMP with scaling of $A$ as follows.

\noindent Let $\lambda_{0}=\mathrm{inf} \{t: K_{X}+B+tA\ \mathrm{is}\ \mathrm{nef} \}$. If $\lambda_{0}<1$, then we stop. Otherwise, by Proposition \ref{findextremalray} there exists a $(K_{X}+ B)$-negative extremal ray $R_{0}$ which satisfies $(K_{X}+B+A)\cdot R_{0}= 0$. We contract this extremal ray. If the contraction
is a Mori fibre contraction, we stop. Otherwise let $\mu_{0}: X\dashrightarrow X_{1}$ be the divisorial contraction or flip. Repeat this process for $(X_{1},\mu_{0\ast}B), \mu_{0\ast}A$ and so on.

\noindent We call this construction a $(K_{X}+B)$-\textit{MMP which is} $(K_{X}+B+A)$-\textit{trivial}.
\edefinition

The following lemma tells us what the output of this construction is if it terminates.

\blemma\label{Output} Let $(X,B)$ be a $\Q$-factorial projective lc pair of dimension $3$ over an algebraically closed field $k$ of characteristic $p>3$, and let $A$ be an effective $\Q$-divisor such that $(X,B+A)$ is lc and $K_{X}+B+A$ is nef. 

\noindent If a $(K_{X}+B)$-MMP which is $(K_{X}+B+A)$-trivial terminates, then its output is a $\Q$-factorial projective lc pair $(X^{\prime}, B^{\prime}+A^{\prime})$, and either

\noindent (1) $X^{\prime}$ is a Mori fibre space over a variety $Y$ , $K_{X^{\prime}}+B^{\prime}+A^{\prime}$ is the pullback of a $\Q$-Cartier $\Q$-divisor from $Y$, and $\mathrm{Supp}\ A^{\prime}$ dominates $Y$, or

\noindent (2) $K_{X^{\prime}}+B^{\prime}+(1-\varepsilon)A^{\prime}$ is nef for any sufficiently small rational $\varepsilon>0$.

\noindent Moreover, $K_{X^{\prime}}+B^{\prime}+A^{\prime}$ is semi-ample if and only if $K_{X}+B+A$ is.
\elemma
\bproof By Definition \ref{MMPKtrivial}, if a $(K_{X}+B)$-MMP which is $(K_{X}+B+A)$-trivial doesn't terminate with a Mori fibre space, then we get a $\Q$-factorial projective lc pair $(X^{\prime}, B^{\prime}+A^{\prime})$ such that $$\lambda:=\mathrm{inf} \{t: K_{X^{\prime}}+B^{\prime}+tA^{\prime}\ \mathrm{is}\ \mathrm{nef} \} <1.$$ 
It is to say that $K_{X^{\prime}}+B^{\prime}+(1-\varepsilon)A^{\prime}$ is nef for any sufficiently small rational $\varepsilon>0$. We are in (2).
If the $(K_{X}+B)$-MMP which is $(K_{X}+B+A)$-trivial terminates with a Mori fibre space $f:(X^{\prime},B^{\prime}+A^{\prime})\to Y$, then by \cite[Theorem 1.4]{hashizume2020minimal} and \cite{hacon2022minimal}, $K_{X^{\prime}}+B^{\prime}+A^{\prime}$ is the pullback of a $\Q$-Cartier $\Q$-divisor from $Y$ since $f$ contracts a $(K_{X^{\prime}}+B^{\prime}+A^{\prime})$-trivial extremal ray.
Moreover, $\mathrm{Supp}\ A^{\prime}$ dominates $Y$ since $f$ only contracts curves which have positive intersections with $A^{\prime}$. We are in (1). The last assertion holds since $K_{X^{\prime}}+B^{\prime}+A^{\prime}$ and $K_{X}+B+A$ coincide after being pulled back to a common resolution of $X'$ and $X$.
\eproof

We will use the following results on termination of flips.

\btheorem\label{Waldrontermination} (\cite[Theorem 1.6]{waldron2018lmmp} and \cite{hacon2022minimal}) Let $(X,B)$ be a projective lc pair of dimension $3$ over a perfect field $k$ of characteristic $p>3$. If $M$ is an effective $\Q$-Cartier $\Q$-divisor on $X$, then any sequence of $(K_{X}+B)$-flips which are also $M$-flips terminates.
\etheorem

\blemma\label{terminaltermination}  Let $(X,B)$ be a $\Q$-factorial projective lc pair of dimension $3$ over an algebraically closed field $k$ of characteristic $p>3$ such that $K_{X}+B$ is nef.
If $X$ is terminal, then any $K_{X}$-MMP which is $(K_{X}+B)$-trivial terminates.
\elemma
\bproof Since every step of a $K_{X}$-MMP which is $(K_{X}+B)$-trivial is a step of a $K_{X}$-MMP, the assertion follows from \cite[Theorem 6.17]{kollar1998birational}.
\eproof

\subsection{Dlt modifications}

The following result helps us to reduce some problems for lc pairs to $\Q$-factorial dlt pairs.

\btheorem\label{dltmodification}  Let $(X, B)$ be an lc pair of dimension $3$ over a perfect field $k$ of characteristic $p>3$.
Then $(X, B)$ has a crepant $\Q$-factorial dlt model. Moreover, we can modify $X$ so that it is terminal.
\etheorem
\bproof The first assertion is proven in the case when $p>5$ in \cite[Theorem 1.6]{birkar2016existence}, and the proof there holds in the case when $p=5$ by \cite{hacon2022minimal}. Let's prove that we can make $X$ terminal. We take a crepant $\Q$-factorial dlt model $g: (X^{\prime},B^{\prime})\to (X,B)$ by the first assertion. Hence, by replacing $(X, B)$ by $(X^{\prime},B^{\prime})$, we may assume that $(X, B)$ is $\Q$-factorial and dlt. Let $U \subseteq X$ be the largest open set such that $(U, B|_{U})$ is a
snc pair. Then $\mathrm{codim}_{ X} (X\backslash U) \geq 2$. 
Note that $(X, 0)$ is klt. By \cite[Theorem 1.7]{birkar2016existence} we take a terminal model
 $f : (X^{\prime}, \Theta^{\prime})\to (X, 0)$  of $(X, 0)$  such that $K_{X^{\prime}} + \Theta^{\prime}= f^{\ast}K_{X}$. Then $f$ is an isomorphism over the smooth locus of $X$; in particular $f$ is an isomorphism
over $U$. Let $Z = X\backslash U$. Define 
$B^{\prime}:= \Theta^{\prime}+ f^{\ast}B$ on $X^{\prime}$
so that
$$K_{X^{\prime}} + B^{\prime} = f^{\ast}(K_{X}+ B),$$
and $(X^{\prime}, B^{\prime})$ is lc.

It remains to show that $(X^{\prime}, B^{\prime})$ is a dlt pair. Let $U^{\prime} = f^{-1}(U)$ and $Z^{\prime} = X^{\prime}\backslash U^{\prime}$.
Then $(U^{\prime}, B^{\prime}|_{U^{\prime}})$ is a snc pair. If $E$ is an exceptional divisor with center in $Z^{\prime}$,
then its center in $X$ is contained in $Z$. Hence $a(E; X^{\prime}, B^{\prime})= a(E; X, B) > -1$.
This completes the proof.
\eproof

\subsection{Some known results on the abundance}

The following theorem collects the recent results towards the abundance conjecture in positive characteristics. 
For (1) see \cite[Theorem 1.3]{waldron2017finite}, \cite[Theorem 3.1]{zhang2019abundance}, \cite[Theorem A]{das2019abundance} and \cite[Theorem 5.1]{xu2024note}). For (2) see \cite[Theorem 5]{witaszek2021canonical} and \cite[Theorem 6.2]{xu2024note}. For (3) see \cite[Theorem 1.1]{zhang2020abundance}, \cite[Corollary
4.13]{witaszek2021canonical} and \cite[Theorem 6.3]{xu2024note}.

\btheorem\label{known} Let $(X,B)$ be a projective lc pair of dimension $3$ over an algebraically closed field $k$ of characteristic $p>3$ such that $K_{X}+B$ is nef. Assume that one of the following conditions holds:

\noindent (1) $\kappa(X,K_{X}+B)\geq 1$,
 
\noindent (2) the nef dimension $n(X,K_{X}+B)\leq 2$,

\noindent (3) the Albanese map $a_{X}: X\to\mathrm{Alb}(X)$ is non-trivial.
 
\noindent Then $K_{X}+ B$ is semi-ample. 
\etheorem

Moreover, the non-vanishing theorem for threefolds in characteristic $p>3$ has been proved. See \cite[Theorem 1.1]{xu2019nonvanishing}, \cite[Theorem 3]{witaszek2021canonical} and \cite[Theorem 4.4]{xu2024note}.

\btheorem\label{nonvan}  Let $(X,B)$ be a projective lc pair of dimension $3$ over a perfect field $k$ of characteristic $p>3$. If $K_{X}+B$ is pseudo-effective, then $\kappa(K_{X}+B)\geq 0$.
\etheorem

\subsection{Slope stability in positive characteristic}

In this subsection, we recall some results on slope stability in positive characteristic.
Let $k$ be an algebraically closed field of characteristic $p>0$.

\bdefinition\label{stabilitydefinition} Let $X$ be a smooth projective variety of dimension $n\geq 2$ over $k$ and let
$D_{1},\cdots,D_{n-1}$ be nef divisors on $X$ such that the $1$-cycle $D_{1}\cdots D_{n-1}$ is not numerically trivial. Let $\sE$ be a torsion free coherent sheaf on $X$. The \textit{slope} of $\sE$ with respect to $(D_{1},\cdots,D_{n-1})$ is defined by 
$$\mu(\sE)=\frac{c_{1}(\sE)\cdot D_{1}\cdots D_{n-1}}{\mathrm{rk}(\sE)}.$$
We call $\sE$ is $(D_{1},\cdots,D_{n-1})$-\textit{semistable} (resp. $(D_{1},\cdots,D_{n-1})$-\textit{stable}) if for every nontrivial subsheaf $\sE^{\prime}\subseteq \sE$ we have
$$\mu(\sE^{\prime})\leq \mu(\sE) (\mathrm{reps.} <).$$
We call $\sE$ is \textit{strongly} $(D_{1},\cdots,D_{n-1})$-\textit{semistable} (resp. \textit{strongly} $(D_{1},\cdots,D_{n-1})$-\textit{stable}) if for every integer $e\geq 0$ the pullback $F^{e\ast}\sE$ is $(D_{1},\cdots,D_{n-1})$-semistable (resp. $(D_{1},\cdots,D_{n-1})$-stable), where $F$ is the Frobenius morphism.
\edefinition

The following theorem is an analogy of Bogomolov inequality in positive characteristic.

\btheorem\label{Bogomolovineqsm} (\cite[Theorem 3.2]{langer2004semistable}) Let $X$ be a smooth projective variety of dimension $n\geq 2$ over $k$ and let
$D_{1},\cdots,D_{n-1}$ be nef divisors on $X$ such that the $1$-cycle $D_{1}\cdots D_{n-1}$ is not numerically trivial. Let $\sE$ be a strongly $(D_{1},\cdots,D_{n-1})$-semistable torsion-free sheaf on $X$. Then
$$\Delta(\sE)\cdot D_{2}\cdots D_{n-1}\geq 0,$$
where $\Delta(\sE):=2\ \mathrm{rk}(\sE) c_{2}(\sE)-(\mathrm{rk}(\sE)-1)c_{1}(\sE)^{2}$.
\etheorem

We will use the following result to prove our main result (see Proposition \ref{specialcase}).

\btheorem \label{constructfoliation} Let $X$ be a smooth projective variety of dimension $n\geq 2$ over $k$ and let
$D_{1},\cdots,D_{n-1}$ be nef divisors on $X$ such that the $1$-cycle $D_{1}\cdots D_{n-1}$ is not numerically trivial. If $c_{1}(T_{X})\cdot D_{1}\cdots D_{n-1}=0$, then $T_{X}$ is strongly $(D_{1},\cdots,D_{n-1})$-semistable if and only if it is $(D_{1},\cdots,D_{n-1})$-semistable. If moreover, $p\geq (n-1)(n-2)$ and $T_{X}$ is not strongly $(D_{1},\cdots,D_{n-1})$-semistable, then the maximal $(D_{1},\cdots,D_{n-1})$-destabiling subsheaf of $T_{X}$ is a $1$-foliation.
\bproof 
The former assertion follows directly from \cite[Theorem 2.1]{mehta2006homogeneous}, once we notice that $\mu_{\mathrm{min}}(T_{X})=\mu(T_{X})=0$ when $T_{X}$ is $(D_{1},\cdots,D_{n-1})$-semistable. The latter assertion follows from the former assertion and \cite[Theorem 3.2]{langer2015generic}.
\eproof
\etheorem

\subsection{$1$-foliations}

We recall the notion of $1$-foliations.

\bdefinition Let $X$ be a smooth variety over an algebraically closed field $k$ of characteristic $p>0$. A $1$-\textit{foliation} on $X$ is a saturated subsheaf $\sF\subseteq T_{X}$ which is involutive (i.e., $[\sF,\sF]\subseteq \sF$) and $p$-closed (i.e., $\xi^{p}\in \sF, \forall \xi\in\sF$). We say that $\sF$ is \textit{smooth} if $\sF$ and $T_{X}/\sF$ are locally free.
\edefinition

The basic properties of $1$-foliations are as follows.

\bproposition (\cite{ekedahl1987foliations}) Let $X$ be a smooth variety over an algebraically closed field $k$ of characteristic $p>0$ and $\sF$ be a $1$-foliation on $X$. Then the following assertions hold.

\noindent (1) $Y:= X/\sF=\mathrm{Spec}_{X}\ Ann(\sF):=\{a\in\sO_{X}|\xi(a)=0, \forall \xi\in\sF\}$ is a normal variety, and there is the following commutative diagram
\beq\begin{tikzcd}
X\arrow[d,swap,"\pi "]\arrow[dr,"F"] \\
Y\arrow[r] & X^{(1)}\notag
\end{tikzcd}
\eeq
where $F$ is the Frobenius morphism, $\pi$ is a morphism with $\mathrm{deg}\ \pi=p^{\mathrm{rk}\ \sF}$.

\noindent (2) There is a one-to-one correspondence between $1$-foliations and normal varieties
between $X$ and $X^{(1)}$. The correspondence is given by $\sF \mapsto X/\sF$ and the inverse correspondence by $Y\mapsto Ann(\sO_{Y}):=\{\xi\in T_{X}|\xi(a)=0, \forall a\in \sO_{Y}\}$.

\noindent (3) $Y$ is regular if and only if $\sF$ is smooth.

\noindent (4) If $Y_{0}$ denotes the regular locus of $Y$ and $X_{0}=\pi^{-1}(Y_{0})$, then 
$$K_{X_{0}}\sim \pi^{\ast}K_{Y_{0}}+(p-1)\mathrm{det}\ \sF|_{X_{0}}.$$
\eproposition

The following criterion helps us to construct $1$-foliations.

\blemma (\cite[Lemma 4.2]{ekedahl1987foliations})\label{foliationcriterion} Let $X$ be a smooth variety over an algebraically closed field $k$ of characteristic $p>0$ and $\sF$ be a saturated $\sO_{X}$-submodule of $T_{X}$. If
$$\mathrm{Hom}(\wedge^{2}\sF,T_{X}/\sF)=\mathrm{Hom}(F^{\ast}\sF,T_{X}/\sF)=0,$$
then $\sF$ is a $1$-foliation.
\elemma

The following theorem tells us how to do ``bend and break" by $1$-foliations.

\btheorem\label{bendbreakfoliation} (\cite[Theorem 2.1]{langer2015generic}) Let $X$ be a smooth variety over an algebraically closed field $k$ of characteristic $p>0$ and $L$ be a nef $\R$-divisor on $X$. Let $f: C \to X$ be a non-constant morphism from a smooth projective curve $C$. Assume that $\sF\subseteq T_{X}$ is a $1$-foliation and smooth along $f(C)$. If
$$c_{1}(\sF)\cdot C > \frac{K_{X}\cdot C}{p-1},$$
then for every $x\in f(C)$ there is a rational curve $B_{x}\subseteq X$ passing through $x$ such that
$$L\cdot B_{x}\leq 2\ \mathrm{dim}\ X\frac{pL\cdot C}{(p-1)c_{1}(\sF)\cdot C-K_{X}\cdot C}.$$

\etheorem

\section{Bertini type theorems in positive characteristic}

In this section, we give some Bertini type results in positive characteristic. We will use these results in Section 5. It is well known that the Bertini theorem fails for basepoint-free linear systems in positive characteristic. However, we have the following special Bertini type result.

\bproposition (\cite[Corollary 4.3 and Proposition 2.4 (1)]{luisa1998axiomatic})\label{specialbertini} Let $f:X\to \P^{n}$ be a finite type morphism from a smooth variety over any infinite field $k$. Assume that $f$ induces separable residue field extensions, i.e. for any point $x\in X$, the residue field of $\P^n$ at $y:=f(x)$ is a separable field extension of the residue field of $X$ at $x$. Then a general section of $|f^{\ast}\sO_{\P^{n}}(1)|$ is smooth.
\eproposition
\bremark Note that we do not require $f$ to be dominant, nor $x\in X$ to be a closed point. The assumption that $f$ induces separable field extension is a very strong condition. For instance, even the projection $g:\A^2_{u,v}\to \A^1_{u}$ does not induce separable residue field extensions: if we take $x$ to be the generic point of the closed subvariety $\mathrm{Spec}(k[u,v]/(u-v^p))$, then the induced residue field extension $k(u)\to k(u^{\frac{1}{p}})$ is not separable. We will apply this proposition in the case where $f$ is a quasi-finite map (see Theorem \ref{bertiniforklt} and Lemma \ref{projectivemodel}). 
\eremark

Next, we follow the strategy in \cite{sato2025general} to prove some Bertini type results in low dimensions.

\bdefinition\label{astcondition} (\cite[Definition 4.1]{sato2025general}) Let $(X,x)$ be a normal surface singularity, and let $f:Y\to X$ be a proper birational morphism from an integral scheme $Y$ with $\mathrm{Exc}(f)=\sum_{i=1}^{n}E_{i}$ such that $f(\mathrm{Exc}(f))=x$. We say that $f$ satisfies the $(\ast)$-\textit{condition} if 

\noindent (1) $E_{i}$ is a smooth curve over the residue field $k(x)$ for every $i$,

\noindent (2)  the scheme theoretic intersection $E_{i}\cap E_{j}$ is smooth over $k(x)$ for every $i\neq j$.
\edefinition

\blemma\label{kltastcondition} Let $(X,x)$ be a normal surface singularity, and let $f:Y\to X$ be the minimal resolution of $(X,x)$  with $\mathrm{Exc}(f)=\sum_{i=1}^{n}E_{i}$ such that $f(\mathrm{Exc}(f))=x$. 
Assume that the following conditions hold:

\noindent (1) the characteristic of the residue field $k(x)$ is greater than $3$,

\noindent (2) $(X,x)$ is a klt singularity,

\noindent (3) $E_{i}$ is geometrically irreducible over $k(x)$ for every $i$.

\noindent Then $f$ satisfies the $(\ast)$-condition.
\bproof Note that $(X,x)$ is an lc and rational singularity (see \cite[Theorem A.3]{sato2025general} for example). Then the assertion follows from the proof of \cite[Proposition 4.4]{sato2025general}.
\eproof
\elemma

\blemma\label{dualgraphunderbasechange} Let $X$ be a normal klt surface over a field $k$ of characteristic $p>3$, and let $x$ be a singular point of $X$. Let $f:Y\to X$ be the minimal resolution of $X$. We write $E_{i}$ for the exceptional prime divisors of $f$ lying over $x$. Assume that the following conditions hold:

\noindent (1) $X$ is smooth over $k$ outside the singular points,

\noindent (2) the residue field $k(x)$ is separable over $k$, 

\noindent (3) $E_{i}$ is geometrically irreducible over $k(x)$ for every $i$.

\noindent Then the surface $X_{\overline{k}}:=X\times_{k}\overline{k}$ is klt at $\tilde{x}$ where $\tilde{x}$ is any point in the preimage of $x$ along the base change. Moreover, the dual graph of $(X_{\overline{k}},\tilde{x})$ is the same as the dual graph of $(X,x)$.
\bproof We may assume that $x$ is the unique singular point of $X$ since $X$ is normal. By (3) and Lemma \ref{kltastcondition}, we have that $f$ satisfies the $(\ast)$-condition. It follows that $E_{i}$ is smooth over $k(x)$ for every $i$. By (2), we know that the residue field of $E_{i}$ at one of its closed points is separable over $k$ for every $i$. Together with (1), it implies that $Y$ is smooth over $k$ 
by \cite[Lemma 00TV]{stacks-project} since $Y$ is regular.
Hence we have the following commutative diagram
\beq\begin{tikzcd}
Y_{\bar{k}}\arrow[d,"\psi"]\arrow[r,"f^{\prime} "]& X_{\bar{k}}\arrow[d,"\phi"]\arrow[r]& \mathrm{Spec}(\bar{k})\arrow[d] \\
Y\arrow[r, "f"] & X\arrow[r]&\mathrm{Spec}(k) \notag
\end{tikzcd}\eeq
such that the following conditions hold:

\noindent  \textcircled{1} $Y_{\bar{k}}= Y\times_{k}\bar{k}$ and $f^{\prime}$ is the base change of $f$ along $\phi$,

\noindent  \textcircled{2} $f^{\prime}$ is a desingularization,

\noindent  \textcircled{3} $K_{Y_{\bar{k}}/X_{\bar{k}}}\sim_{\Q}\psi^{\ast}K_{Y/X}$,

\noindent where \textcircled{2} follows from the smoothness of $Y$, \textcircled{3}  follows by \cite[Lemma 0AWD]{stacks-project}.

Since $f$ is the minimal resolution of $X$, $K_{Y/X}$ is nef. Hence by \textcircled{3}, $K_{Y_{\bar{k}}/X_{\bar{k}}}$ is nef. Together with \textcircled{2}, it implies that $f^{\prime}$ is the minimal resolution of $X_{\bar{k}}$. We write $\tilde{x}$ for a point in the preimage of $x$ along $\phi$ and $F_{i}\subseteq \psi^{-1}(E_{i})$ be the exceptional prime divisors over $\tilde{x}$. Then by Definition \ref{dualgraphdefinition} the dual graph of $(X_{\bar{k}},\tilde{x})$ is determined by $\mathrm{dim}_{k(\tilde{x})}H^{0}(F_{i},\sO_{F_{i}})$, $g(F_{i})$ and $F_{i}\cdot F_{j}$. 
Since  $E_{i}$ are smooth over $k(x)$,
the coefficients of $\psi^{\ast}E_{i}$ at $F_{i}$ are one.
Since the base change 
$$\mathrm{Spec}(k(\tilde{x}))=\mathrm{Spec}(\bar{k})\to \mathrm{Spec}(k(x))$$ 
is flat and 
$$F_{i}=E_{i}\times_{k(x)}k(\tilde{x}),$$ 
we have
$$\mathrm{dim}_{k(\tilde{x})}H^{0}(F_{i},\sO_{F_{i}})=\mathrm{dim}_{k(x)}H^{0}(E_{i},\sO_{E_{i}})$$
and by Definition \ref{genusdefinition}
$$g(F_{i})=\frac{\mathrm{dim}_{k(\tilde{x})}H^{1}(F_{i},\sO_{F_{i}})}{\mathrm{dim}_{k(\tilde{x})}H^{0}(F_{i},\sO_{F_{i}})}=\frac{\mathrm{dim}_{k(x)}H^{1}(E_{i},\sO_{E_{i}})}{\mathrm{dim}_{k(x)}H^{0}(E_{i},\sO_{E_{i}})}=g(E_{i}).$$ 
Moreover, by Definition \ref{intersectionnumber} 
$$F_{i}\cdot F_{j}= F_{i}\cdot \psi^{\ast}E_{j}=\mathrm{deg}_{F_{i}/k(\tilde{x})}(\psi^{\ast}E_{j})=\mathrm{deg}_{E_{i}/k(x)}(E_{j}) = E_{i}\cdot E_{j}.$$
Therefore, the assertion holds.
\eproof
\elemma

The following theorem is a Bertini type result which tells us what singularities general hyperplane sections of klt quasi-projective threefolds have in characteristic $p>3$. See also \cite{sato2025general,sato2020general}.

\btheorem\label{bertiniforklt} Let $X$ be a klt quasi-projective threefold over an algebraically closed field $k$ of characteristic $p>3$. Let $C$ be a curve contained in the singular locus of $X$, and let $H$ be a general hyperplane section of $X$.  Then, for a closed point $x\in C\cap H$,  we have $(H,x)$ is klt. Moreover, there exists a finite cover $g:V\to X$ such that $g$ is \'{e}tale over general points in $C$, and the dual graph of $(H,x)$ is the same as the dual graph of $(V_{\gamma},\gamma)$,  where $\gamma$ is the generic point of a curve $\Gamma\subseteq V$ lying over $C$, and $V_{\gamma}$ is the localization of $V$ at $\gamma$.
\bproof Note that once we have a finite \'{e}tale cover $V_U$ of an open subset $U\subseteq X$ containing the generic point of $C$, we may take $V$ to be the normalization of $X$ in the functional field $k(V_U)$.  Hence the problem is local, and we can always shrink $X$ around the generic point of $C$ if necessary. 
We assume that $X$ is affine and the singular locus of $X$ is $C$. Let $\delta$ be the generic point of $C$, and let $f: Y_{\delta}\to X_{\delta}$ be the minimal resolution of $(X_{\delta},\delta)$, where $\mathrm{Exc}(f)=\sum_{i}F_{i}$ and $X_{\delta}$ is the localization of $X$ at $\delta$. Let $\pi:X\to T$ be a  projection induced by a general pencil of hyperplane sections, where $T$ is an open set of $\P^{1}$.
Then, $C$ dominates $T$, and
$\delta=\mathrm{Spec}(k(C))$ is mapped to $\eta:=\mathrm{Spec}(k(T))$.
Let $L_{1}/k(C)$ be a finite separable field extension such that every irreducible component of $F_{i}\times_{k(C)}L_{1}$ is geometrically irreducible over $L_{1}$ for every $i$.  After shrinking $T$, the normalization $\phi:T^{\prime}\to T$ of $T$ in $L_{1}$ is an \'{e}tale map. Let 
$$\pi^{\prime}:V\to T^{\prime}$$ 
be the base change of $\pi$ along $\phi$. Let $\Gamma\subseteq V$ be a curve lying over $C$, and let $\gamma$ denote the generic point of $\Gamma$. Then we have that $\mathrm{Spec}(k(\Gamma))$ is a point in
$$\mathrm{Spec}\big( k(C)\otimes_{k(T)}L_{1}\big).$$
Since $k(C)\otimes_{k(T)}L_{1}$ is a vector space over $L_{1}$, it follows that 
$$L_{1}\subseteq k(\Gamma).$$
Note that $V$ is finite over $X$. Hence up to replacing $X$ by $V$, $T$ by $T'$, $\pi$ by $\pi^{\prime}$, $C$ by $\Gamma$ and $\delta$ by $\gamma$, we can assume that 
every exceptional prime divisor $F_{i}$ of the minimal resolution $f$ is geometrically irreducible over $k(C)$. Moreover, to prove the assertion, it suffices to show that the dual graph of $(H,x)$ is the same as the dual graph of $(X_{\delta},\delta)$.

Let $X_{\eta}$ be the generic fibre of $\pi$.  By  \cite[Theorem 5.2]{sato2025general}, $X_{\eta}$ is geometrically klt over $k(T)$. Note that $\delta$ is a closed point in $X_{\eta}$.
 Moreover, the natural map
$$X_{\delta}\to  X_{\eta}$$
induced by the universal property of $X_{\eta}$ is isomorphic at $\delta$. It implies that the dual graph of $(X_{\delta},\delta)$ is the same as the dual graph of $(X_{\eta},\delta)$.
Since $X\backslash C$ is smooth over $k$, by Proposition \ref{specialbertini} general fibres of $\pi:X\backslash C\to T$ are smooth over $k$. By \cite[Proposition 2.4 (1)]{luisa1998axiomatic} it follows that $X_{\eta}\backslash \delta$ is smooth over $k(T)$. Since $\pi$ is induced by a general pencil of hyperplane sections and an \'{e}tale base change, the finite field extension $k(C)/k(T)$ is separable. Moreover, every exceptional prime divisor $F_{i}$ of the minimal resolution $f$ is geometrically irreducible over $k(C)$ by our assumption. Hence by Lemma \ref{dualgraphunderbasechange} the dual graph of $(X_{\eta},\delta)$ is the same as the dual graph of $(X_{\overline{\eta}},\tilde{\delta})$, where $X_{\overline{\eta}}:=X_{\eta}\times_{k(T)}\overline{k(T)}$ and $\tilde{\delta}$ is any point in the preimage of $\delta$ along the base change. It implies that the dual graph of $(X_{\delta},\delta)$ is the same as the dual graph of $(X_{\overline{\eta}},\tilde{\delta})$.

We prove that the dual graph of $(X_{\overline{\eta}},\tilde{\delta})$ is the same as the dual graph of $(H,x)$. After shrinking $X_{\overline{\eta}}$, we can assume that $\tilde{\delta}$ is the unique closed point in $X_{\overline{\eta}}$ lying over $\delta$.
Let $\rho_{\overline{\eta}}: Y_{\overline{\eta}} \to X_{\overline{\eta}}$ be the minimal resolution of $X_{\overline{\eta}}$. Then there exist a finite field extension $L_{2}$ of $k(T)$ and a Cartesian diagram
\beq\begin{tikzcd}
Y_{\overline{\eta}}\arrow[d]\arrow[r,"\rho_{\overline{\eta}} "]& X_{\overline{\eta}}\arrow[d] \\
Y_{L_{2}}\arrow[r, "\rho_{L_{2}}"] & X_{L_{2}}:=X_{\eta}\times_{k(T)}L_{2} \notag
\end{tikzcd}\eeq
where $\rho_{L_{2}}$ is the minimal resolution of $X_{L_{2}}$ and horizontal arrows are base changes of $f$ along the base field extensions. Moreover, by Lemma \ref{kltastcondition} we can assume that every stratum of 
$$\mathrm{Exc}(\rho_{L_{2}})=\sum_{i=1}^{n}E_{i,L_{2}}$$ 
is smooth over $L$. Then by Definition \ref{dualgraphdefinition} the dual graph of $(X_{\overline{\eta}},\tilde{\delta})$ is determined by $E_{i,L_{2}}\cdot E_{j,L_{2}}$.
We replace $k(T)$ by $L_{2}$, $T$ by its normalization in $L_{2}$ and $X$ by its base change along the normalization. Since $H$ is a fibre of $\pi$, $(H,x)$ doesn't change under this replacement.
Then by taking integral models there exists a following commutative diagram
\beq\begin{tikzcd}
Y_{\eta}\arrow[d]\arrow[r,"\rho_{\eta} "]& X_{\eta}\arrow[d] \\
Y_{U_{1}}\arrow[r, "\rho_{U_{1}}"] & X_{U_{1}} \notag
\end{tikzcd}\eeq
where $U_{1}\subseteq T$ is an open and dense subset, $X_{U_{1}}:=\pi^{-1}(U_{1})$, $Y_{U_{1}}$ is a quasi-projective threefold such that $\rho_{U_{1}}$ is a projective morphism, and vertical arrows are natural injections. After shrinking $X$ and $T$, we can assume that $U_{1}=T$ and $Y:=Y_{T}$ is smooth over $T$, since $Y$ is smooth over $\eta$.
Moreover, by the following claim, we can assume that $\rho_{t}:Y_{t}\to X_{t}$ is a minimal resolution for every $t\in T$, where $Y_{t},X_{t}$ are fibres of $Y,X$ over $t$.

\bclaim There exists an open and dense subset $U_{2}\subseteq T$ such that $K_{(\pi\circ\rho_{T})^{-1}(U_{2})}$ is nef over $\pi^{-1}(U_{2})$.
\eclaim
\noindent \textit{Proof of the claim.} By Theorem \ref{lcmmp}, we can run a $K_{Y}$-MMP over $X$ to get a relative minimal model $Y^{\prime}$. Then $K_{Y^{\prime}/X}$ is nef. Since $Y^{\prime}$ is terminal, it has isolated singularities. Hence after shrinking $X$ and $T$, we can assume that $Y^{\prime}$ is smooth over $T$. Note that the natural map $Y^{\prime}\to X$ gives the minimal resolution of $X_{\eta}$. Hence it is isomorphic to $\rho_{\eta}$ over $\eta$ since  $\rho_{\eta}:Y_{\eta}\to X_{\eta}$ is the minimal resolution. It implies that every step the MMP only contract subvarieties which don't dominate $T$. Hence, the assertion holds.  $\hfill\square$\\

We write 
$$\mathrm{Exc}(\rho_{T})=\sum_{i=1}^{n}E_{i}$$ 
such that $E_{i}|_{Y_{\eta}}=E_{i,k(T)}$. 
After shrinking $X$ and $T$, we can assume that $E_{i}$ is smooth over $T$ for every $i$ since $E_{i,k(T)}$ is smooth over $\eta$.
Then $E_{i,t}$ are all prime $\rho_{t}$-exceptional divisors for all $t\in T$, where $E_{i,t}$ are the restrictions of $E_{i}$ on the fibres $Y_{t}$.
It implies that 
$$E_{i,t}\cdot E_{j,t}=E_{i}|_{Y_{t}}\cdot E_{j}|_{Y_{t}}=(E_{i}\cdot E_{j})|_{Y_{t}}=(E_{i}\cdot E_{j})|_{Y_{\eta}}=E_{i}|_{Y_{\eta}}\cdot E_{j}|_{Y_{\eta}}=E_{i,k(T)}\cdot E_{j,k(T)}$$
 for all $t\in T$, where the third equality holds since both $(E_{i}\cdot E_{j})|_{Y_{t}}$ and $(E_{i}\cdot E_{j})|_{Y_{\eta}}$ equal to the degree of $E_{i}\cdot E_{j}$ over $T$ as a $1$-cycle.
Hence the dual graph of $(X_{\overline{\eta}},\tilde{\delta})$ is the same as the dual graph of $(H,x)$.
In conclusion, the assertion holds.
\eproof
\etheorem
\bcorollary\label{bertiniforcyclicquotient} Let $(X,B)$ be an lc quasi-projective threefold over an algebraically closed field $k$ of characteristic $p>3$ with $B$ is reduced. Let $C\subseteq B$ be a curve contained in the singular locus of $X$ and $H$ be a general hyperplane section of $X$.  
Assume that $(X,B)$ is plt at $\delta$ and $(X_{\delta},\delta)$ is a cyclic quotient singularity whose index is $m$,  where $\delta$ is the generic point of $C$ and $X_{\delta}$ is the localization of $X$ at $\delta$.
Then for a closed point $x\in C\cap H$,  we have $(H,x)$ is a cyclic quotient singularity whose index is $m$.
\bproof Since $\delta \in B\neq 0$, by Proposition \ref{lcsurfacesing}, every exceptional prime divisor of the minimal resolution of $(X_{\delta},\delta)$ is isomorphic to $\P^{1}_{k(C)}$. Hence it is geometrically irreducible over $k(C)$. It implies that the dual graph of $(X_{\delta},\delta)$ doesn't change under any \'{e}tale base change. Moreover, since $X$ is klt along $C$, we can apply Theorem \ref{bertiniforklt} to conclude that the dual graph of $(X_{\delta},\delta)$ is as same as the dual graph of $(H,x)$. Hence, the assertion holds.
\eproof
\ecorollary

Since we don't know whether the Bertini theorem holds for the pullback of a very ample linear system along a birational morphism, we prove the following special result. We will use it in the proof of Proposition \ref{keyinequality}.

\bproposition\label{gooddesingularization} Let $X$ be a klt quasi-projective threefold over an algebraically closed field $k$ of characteristic $p>3$. Then there exists a desingularization $\rho:V\to X$ such that for a general hyperplane section $H$ of $X$, we have that $S:=\rho^{-1}(H)$ is smooth, and $\rho|_{S}:S\to H$ is the minimal resolution of $H$.
\bproof Let $\rho_{1}: V_{1}\to X$ be a desingularization such that $\rho_{1}$ is an isomorphism outside the singular locus of $X$.  By Theorem \ref{lcmmp}, we can run a $K_{V_{1}}$-MMP over $X$ to get a relative minimal model $V_{2}$. Then $V_{2}$ is terminal and $K_{V_{2}}$ is nef over $X$. 
There is a morphism $\rho_{2}:V_{2}\to X$.
Let $\rho_{3}:V\to V_{2}$ be a desingularization such that $\rho_{3}$ is an isomorphism outside the singular locus of $V_{2}$. We prove that 
$$\rho: V\overset{\rho_{3}}\to V_{2}\overset{\rho_{2}}\to X$$ satisfies the condition of the assertion. Let $H$ be a general hyperplane section of $X$. Then by Theorem \ref{bertiniforklt}, $H$ has klt singularities. Let $C$ be a curve contained in the singular locus of $X$ and $h\in H\cap C$ be a closed point. Then for a sufficiently small open neighbourhood $h\in U\subseteq X$, we have 
$$\rho^{-1}(U)=\rho_{2}^{-1}(U)$$ is smooth since $V_{2}$ is terminal, $H$ is general and $\rho_{3}$ is an isomorphism outside the singular locus of $V_{2}$. Hence it suffices to show that 
$$S:=\rho_{2}^{-1}(U\cap H)$$ is smooth. 
This problem is local on $U$. After shrinking $U$, by Theorem \ref{bertiniforklt} we can replace $U$ by  its \'{e}tale base change so that the dual graph of $(H,h)$ is the same as the dual graph of $(U_{\delta},\delta)$,  where $\delta$ is the generic point of $C$, and $U_{\delta}$ is the localization of $U$ at $\delta$.
We replace $X$ by $U$, $H$ by $H\cap U$, $V_{2}$ by $\rho^{-1}_{2}(U)$ and define $\rho^{\prime}:=\rho_{2}|_{S}$.

Shrinking $X$ if necessary, we assume $H\sim_{\Q}0$. Thus $S=\rho_{2}^{\ast}H \sim_{\Q} 0$ on $V_{2}$. It follows that
$$K_{V_{2}}|_{S}\sim_{\Q}(K_{V_{2}}+S)|_{S}\sim_{\Q}K_{S}.$$
Hence $K_{S}$ is nef over $H$ since $K_{V_{2}}$ is nef over $X$. 
We write 
$$K_{S}\sim_{\Q}\rho^{\prime\ast}K_{H}+\sum_{i}a_{i}E_{i}.$$
Then by the negativity lemma \cite[Lemma 3.39]{kollar1998birational}, we have 
$$a_{i}\leq 0.$$
It is to say that $\rho^{\prime}$ doesn't contract any exceptional prime divisor $E$ over $H$ with the discrepancy $a(E;H)>0$. To prove that $S$ is smooth, i.e. terminal, it suffices to show that $\rho^{\prime}$ contracts all exceptional prime divisors $E$ over $H$ with the discrepancies $a(E;H)\leq 0$.
Note that $\rho_{2}$ restricts to the minimal resolution 
$$\rho_{\delta}: \rho_{2}^{-1}(X_{\delta})\to X_{\delta}$$ 
of $X_{\delta}$ with the exceptional divisor $\mathrm{Exc}(\rho_{\delta})=\Sigma_{j=1}^{r}F_{\delta,j}$. 
For every $F_{\delta,j}$, there exists an exceptional prime divisor $F_{j}$ of $\rho_{2}$ which restricts to $F_{\delta,j}$.
It follows that $\rho_{2}$, and hence $\rho^{\prime}$ contract at least $r$ exceptional prime divisors. By our assumption, the dual graph of $(H,h)$ is the same as the dual graph of $(X_{\delta},\delta)$. It implies that $\rho^{\prime}$ contracts all exceptional prime divisors $E$ over $H$ with the discrepancies $a(E;H)\leq 0$.
Therefore, the assertion holds.
\eproof
\eproposition

\section{Chern classes on normal varieties in positive characteristic}

In this section, we recall Langer's definition of the Chern classes of reflexive sheaves on normal varieties in positive characteristic (see Subsection \ref{subsection 5.2} and \ref{subsection 5.3}). Moreover, using Lefschetz-Riemann-Roch theorem, we calculate local contributions of some quotient singularities to Riemann-Roch formula for normal surfaces. Then we give lower bounds of some relative second Chern classes (see Subsection \ref{subsection 5.4}).

\subsection{Local relative Chern classes}\label{subsection 5.2}

In this subsection, we recall the notion of local relative Chern classes for resolutions of normal surfaces. The main reference of this subsection is \cite[Section 3]{langer2025intersection}. 
Let $k$ be an algebraically closed field and let $A$ be an excellent normal $2$-dimensional Henselian local $k$-algebra. Let $X =\mathrm{Spec}\ A$ and let $x\in X$ be the closed point of $X$. Let $f: \tilde{X}\to X$ be a desingularization of $X$ with $\mathrm{Exc}(f)=\Sigma_{i}E_{i}$.

\bdefinition\label{relc1definition} Let $\sF$ be a vector bundle on $\tilde{X}$. We define the \textit{first relative Chern class} of $\sF$ with respect to $f$ to be the unique $\Q$-divisor $c_{1}(f,\sF)$ supported on $E$ such that for every irreducible component $E_{i}$ of $E$ we have
$$c_{1}(f,\sF)\cdot E_{i}=\mathrm{deg}\ \sF|_{E_{i}}.$$
\edefinition
\bremark The existence and uniqueness of $c_{1}(f,\sF)$ follow from the fact that the intersection matrix $[E_{i}\cdot E_{j}]$ of the exceptional divisors is negative definite (see \cite[3.1]{langer2025intersection}).
\eremark

Let $\tilde{\pi}: \tilde{Y}\to \tilde{X}$  be a generically finite proper morphism from a regular surface $\tilde{Y}$. The Stein factorization of $f\circ \tilde{\pi}$ gives a proper birational morphism $g:\tilde{Y}\to Y$. Then we have the following commutative diagram
\beq\begin{tikzcd}
\tilde{Y}\arrow[d,"\tilde{\pi}"]\arrow[r,"g "]& Y\arrow[d,"\pi"] \\
\tilde{X}\arrow[r, "f"] & X. \notag
\end{tikzcd}\eeq
\noindent Fix a rank $r$ vector bundle $\sF$ on $\tilde{X}$. By \cite[Theorem 4.7 (ii)]{kleiman1969geometry} possibly after further blowing up $\tilde{Y}$, there exists a filtration
$$ 0=\sF_{0}\subseteq \sF_{1} \cdots\subseteq \sF_{r}=\tilde{\pi}^{\ast}\sF $$
 such that all quotients $\sL_{i} = \sF_{i}/\sF_{i-1}$ are line bundles on $\tilde{Y}$. We define $L_{i}:=c_{1}(g,\sL_{i})$.

\bdefinition\label{relc2definition} The \textit{second relative Chern class} $c_{2}(f,\sF)$ of $\sF$ with respect to $f$
is defined as the real number
$$c_{2}(f,\sF):=\mathrm{inf}\Big(\frac{\sum_{i<j}L_{i}L_{j}}{\mathrm{deg}\ \tilde{\pi}} \Big),$$
where the infimum is taken over all $\tilde{\pi}$ and filtrations as above.
\edefinition

\bdefinition\label{relative_Euler_characteristic} We define the \textit{relative Euler characteristic} $\chi(f,\sF)$ by
$$\chi(f,\sF):=\mathrm{dim}\ H^{0}(X,f_{[\ast]}\sF/f_{\ast}\sF)+\mathrm{dim}\ H^{0}(X,R^{1}f_{\ast}\sF)$$
and we set
$$a(f,\sF):=\chi(f,\sF)-r \chi(f,\sO_{\tilde{X}})+\frac{1}{2}c_{1}(f,\sF)(c_{1}(f,\sF)-K_{\tilde{X}})-c_{2}(f,\sF).$$
\edefinition
\bremark\label{local_contribution_definition} The invariant $a(f,\sF)$ is important since it doesn't depend on the choice of desingularization $f$ (see \cite[Proposition 3.7]{langer2025intersection}). Hence, for a reflexive sheaf $\sE$ on $X$,  we can set
$$a(x,\sE):=a(f',\sF'),$$
where $f': X' \to X$ is any desingularization of $X$ and $\sF'$ is any vector bundle on $X'$ such that 
$f'_{[\ast]}\sF'\cong \sE$.
\eremark
 
We will use the following result in the proof of Proposition \ref{keyinequality}.

\bproposition\label{sumprinciple} Assume that the base field $k$ has characteristic $p>0$, and that $\sF_{1}$, $\sF_{2}$, and $\sF_{3}$ are vector bundles on $\tilde{X}$.
If there exists a  following exact sequence
$$0\to \sF_{1}\to\sF_{2}\to\sF_{3}\to 0,$$
then 
$$\mathrm{ch}_{2}(f,\sF_{2})=\mathrm{ch}_{2}(f,\sF_{1})+\mathrm{ch}_{2}(f,\sF_{3})$$
where $\mathrm{ch}_{2}(f,\cdot):=\frac{1}{2}c_{1}(f,\cdot)^{2}-c_{2}(f,\cdot)$.
\bproof By \cite[Theorem 3.15]{langer2025intersection}, for $i=1,2,3$ and any positive integer $m$ we have
$$\mathrm{ch}_{2}(f,\sF_{i})=\frac{1}{2}c_{1}(f,\sF_{i})^{2}-c_{2}(f,\sF_{i})=-\mathrm{lim}_{m\to\infty}\frac{\chi(f,(F^{m})^{\ast}\sF_{i})}{p^{2m}}.$$
We consider the exact sequence $0\to \sF_{1}\to\sF_{2}\to\sF_{3}\to 0$.
Since $\tilde{X}$ is smooth, the maps $F^m$ are flat.
Hence we have the exact sequences
$$0\to(F^{m})^{\ast} \sF_{1}\to(F^{m})^{\ast}\sF_{2}\to(F^{m})^{\ast}\sF_{3}\to 0.$$
Thus the assertion follows from the additivity of the relative Euler characteristic.
\eproof
\eproposition

\subsection{Chern classes of reflexive sheaves}\label{subsection 5.3}

In this subsection, we recall the definition of Chern classes of reflexive sheaves on normal varieties. The main reference of this subsection is \cite[Section 4 and 5]{langer2025intersection}.
Let $X$ be a proper normal surface defined over an algebraically closed field $k$, and let $\sE$ be a reflexive sheaf on $X$.  
Let $f :\tilde{X} \to X$ be a desingularization with $\mathrm{Exc}(f)= E$.
For every $x\in f(E)$ we consider the map $\nu_{x}:\mathrm{Spec}\ \sO^{h}_{X,x} \to X$ from the spectrum of the henselization of the local ring of $X$ at $x$ and the base change 
\beq f_{x} = \nu^{\ast}_{x}f : \tilde{X}_{x} \to \mathrm{Spec}\ \sO^{h}_{X,x} \eeq
of $f$ via $\nu_{x}$. Note that the group of divisors on $\tilde{X}_{x}$ that are supported on
the exceptional locus $E_{x}$ of $f_{x}$ embeds into $A_{1}(\tilde{X})$, where $A_{1}(\tilde{X})$ is the group of $1$-cycles modulo rational equivalence on $\tilde{X}$. For any vector bundle $\sF$ on $\tilde{X}$ we write $c_{1}(f_{x},\sF)\in A_{1}(\tilde{X})\otimes \Q$ for the image of $c_{1}(f_{x},\sF|_{\tilde{X}_{x}})$. Then we have
$$c_{1}(\sF)=f^{\ast}c_{1}(f_{[\ast]}\sF)+\sum_{x\in f(E)}c_{1}(f_{x},\sF)\in A_{1}(\tilde{X})\otimes\Q.$$
We write $a(x,\sE)$ for $a(x,\nu^{\ast}_{x}\sE)$ (see Remark \ref{local_contribution_definition}).

\bdefinition\label{c2definition} Let $\sF$ be a vector bundle on $\tilde X$ such that $f_{[\ast]}\sF\cong \sE$. Using the homomorphism $f_{\ast}:A_{0}(\tilde{X})\to A_{0}(X)$, we define the \textit{second Chern class} of $\sE$ to be 
$$c_{2}(\sE):=f_{\ast}c_{2}(\sF)-\sum_{x\in f(E)}c_{2}(f_{x},\sF)[x]\in A_{0}(X)\otimes\R,$$
where $A_{0}(\tilde{X})$ and $A_{0}(X)$ denote the groups of $0$-cycles modulo rational equivalence on $\tilde{X}$ and $X$, respectively.
Moreover, we define the \textit{discriminant}  
$$\Delta(\sE):=2r\ c_{2}(\sE)-(r-1)c_{1}(\sE)^{2}\in A_{0}(X)\otimes\R$$
where $r$ is the rank of $\sE$. For a more detailed explanation of these definitions,  see \cite[4.1]{langer2025intersection}.
\edefinition
 
In higher dimensions, we can define Chern classes by reducing to the case when $X$ is a surface (see \cite[Section 5]{langer2025intersection}). 
We will use the following result in Proposition \ref{specialcase}.

\btheorem (\cite[Theorem 3.4]{langer2024bogomolov})\label{Bogomolovineq} 
Let $X$ be a normal projective threefold over an algebraically closed field of characteristic $p>0$.
Let
$L_{1}$ and $ L_{2}$ be nef line bundles on $X$ such that $L_{1}\cdot L_{2}$ is numerically nontrivial, and  
let $\sE$ be a strongly $(L_{1},L_{2})$-semistable reflexive sheaf on $X$. Then
$$\Delta(\sE)\cdot L_{2}:=2\mathrm{rk}(\sE)\int_{X} c_{2}(\sE)\cdot L_{2}-(\mathrm{rk}(\sE)-1)\int_{X} c_{1}(\sE)^{2}\cdot L_{2}\geq 0,$$
where $\int_{X}$ is the degree map of $0$-cycles on $X$ (see \cite[Definition 1.4]{fulton2013intersection}).
\etheorem

\subsection{Computations of some relative second Chern classes}\label{subsection 5.4}

In this subsection, we compute relative second Chern classes of some cyclic quotient singularities.
\bsetting\label{setting}
Let $k$ be an algebraically closed field of characteristic $p>0$, and let $X =\hat{\mathbb{A}}_{k}^{2}/\Z_{m}$ (see Remark \ref{cyclic_quotient_action}), where $m$ is a positive integer not divisible by $p$. Let $f: \tilde{X}\to X$ be a desingularization of $X$ with $\mathrm{Exc}(f)=\Sigma_{i}E_{i}$ such that $\mathrm{Exc}(f)$ is snc. We write $\Omega_{X}^{[1]}:=(\Omega_{X})^{\ast\ast}$. We will compute $a(0,\Omega_{X}^{[1]})$ (see Remark \ref{local_contribution_definition}), where $0$ is the closed point of $X$, to give lower bounds of $c_{2}(f,\Omega_{\tilde{X}})$.
\esetting
\blemma\label{projectivemodel} There is a smooth projective variety $T$ over $k$ and a group monomorphism $\Z_{m}\emb \mathrm{Aut}\ T$ such that the following assertions hold.

\noindent (1) $\big|\mathrm{Fix}(\Z_{m})\big|< +\infty$, where $\mathrm{Fix}(\Z_{m})$ is the set of fixed points for the action of $\Z_{m}$.

\noindent (2) For any $t\in \mathrm{Fix}(\Z_{m})$ we have $\hat{T}_{t}/\Z_{m}\cong X$, where $\hat{T}_{t}$ is the completion of $T$ at $t$.
\bproof The proof is similar to \cite[P.405]{reid1985young}. By the definition, $X=\hat{\mathbb{A}}_{k}^{2}/\Z_{m}$ comes from a linear and diagonalized action of $\Z_{m}$ on $\A^{2}$. This action gives an action of $\Z_{m}$ on $\A^{5}$ by 
$$g\in \Z_{m}:(x_{1},x_{2})\in \A^{2}\times \A^{3}\cong\A^{5} \mapsto (g(x_{1}),x_{2}).$$
It induces a natural action of $\Z_{m}$ on $\P^{4}$. Then it's clear that the fixed point set 
$$\mathrm{Fix}(\Z_{m})\subseteq (x_{1}=0)\cup(x_{2}=0)\subset\P^{4}.$$
Consider the quotient map $\pi:\P^{4}\to\P^{4}/\Z_{m}$. We use two general hyperplane sections to cut out a subvariety $Y\subseteq\P^{4}/\Z_{m}$. By Proposition \ref{specialbertini}, $T:=\pi^{-1}(Y)$ is smooth. Since $\mathrm{dim}(\mathrm{Fix}(\Z_{m}))=2$, we have (1) holds. Moreover, (2) holds by our construction.
\eproof
\elemma

In order to compute $a(0,\Omega_{X}^{[1]})$, we need to apply the Lefschetz fixed point formula in positive characteristic.  Let $W(k)$ be the Witt ring of $k$.
For $c\in k$, we set $w(c):=(c,0,0,\cdots)\in W(k)$. 
Let $\psi$ be an endomorphism of a finite rank $k$-vector space. Recall that, the \textit{Brauer trace} $tr_{B}(\psi)$ of $\psi$, is defined to be $\sum_{i} w(c_{i})$,
where $c_{i}$ are eigenvalues of $\psi$, with the appropriate multiplicities. We refer to \cite[18.1]{serre1977linear} for more details.

\blemma\label{fixedpointformula} Let $T$ be a smooth projective surface over $k$ with an action of $\Z_{m}$ on $T$. Assume that $\big|\mathrm{Fix}(\Z_{m})\big|< +\infty$, where $\mathrm{Fix}(\Z_{m})$ is the set of fixed points for the action of $\Z_{m}$. For any $\Z_{m}$-equivariant locally free sheaf $\sF$ on $T$ and $0\neq g\in \Z_{m}$, we have
$$\sum_{i}(-1)^{i}tr_{B}(g^{\ast}|_{H^{i}(T,\sF)})=\sum_{t\in\mathrm{Fix}(\Z_{m})}\ \frac{tr_{B}(g|_{\sF_{t}})}{(1-w(c_{g,1}))(1-w(c_{g,2}))},$$
where  $\sF_{t}$ are the fibre of $\sF$ at $t$ and $c_{g,i}$ are the eigenvalues of $g$ on $\A^{2}$.
\bproof It follows from the Lefschetz fixed point formula in positive characteristic (see \cite[0.6]{baum1979lefschetz} for example).
\eproof
\elemma

\blemma\label{Brauertr} (\cite[18.1]{serre1977linear}) Let $V$ be a vector space over $k$ and let $\Z_{m}\to \mathrm{GL}(V)$ be a representation. Then we have
$$\mathrm{dim}(V^{\Z_{m}})=\frac{1}{m}\sum_{g\in\Z_{m}}tr_{B}(g),$$
where $V^{\Z_{m}}$ is the subspace of $V$, which is fixed by $\Z_{m}$.
\elemma

The following result is an adaptation of \cite[Theorem 5.4]{langer2000chern}.

\btheorem\label{localcontribution} Let $X$ be as stated in Setting \ref{setting}. Then we have 
$$a(0,\Omega_{X}^{[1]})=\frac{1}{m}-1.$$
\bproof 
Let $T$ be a smooth projective surface constructed in Lemma \ref{projectivemodel}, and let $\pi:T\to Y:=T/\Z_{m}$ be the quotient map.
We set $N:=|\mathrm{Fix}(\Z_{m})|$ and $y:=\pi(t)$, where $t\in \mathrm{Fix}(\Z_{m})$ is a fixed point for the action of $\Z_{m}$. Note that by (2) of Lemma \ref{projectivemodel}, $a(y,\Omega_{Y}^{[1]})=a(0,\Omega_{X}^{[1]})$ is independent of the choice of $t$.
It follows from (4) of \cite[Proposition 4.2]{langer2025intersection} that $\int_{Y}c_{2}(\Omega_{Y}^{[1]})=\frac{1}{m}\int_{Y}c_{2}(\Omega_{T})$.
Moreover, by \cite[Theorem 4.4]{langer2025intersection}, we have
$$\chi(Y,\Omega_{Y}^{[1]})=-\int_{Y}c_{2}(\Omega_{Y}^{[1]})+2\chi(Y,\sO_{Y})+N\ a(y,\Omega_{Y}^{[1]}),$$
and
$$\int_{Y}c_{2}(\Omega_{T})=-\chi(T,\Omega_{T})+2\ \chi(T,\sO_{T}).$$
It implies that
\beq\begin{split}
a(y,\Omega_{Y}^{[1]})=&\frac{1}{N}\big( \chi(Y,\Omega_{Y}^{[1]})+\int_{Y}c_{2}(\Omega_{Y}^{[1]})-2\ \chi(Y,\sO_{Y}) \big)\\
=&\frac{1}{N}\big(\sum_{i}(-1)^{i}\mathrm{dim}(H^{i}(Y,\Omega_{Y}^{[1]}))+\frac{1}{m}\int_{Y}c_{2}(\Omega_{T})-2 \sum_{i}(-1)^{i}\mathrm{dim}(H^{i}(Y,\sO_{Y})) \big)\\
=&\frac{1}{N}\big(\sum_{i}(-1)^{i}\mathrm{dim}(H^{i}(T,\Omega_{T})^{\Z_{m}})+\frac{1}{m}(-\chi(T,\Omega_{T})+2\ \chi(T,\sO_{T}))\\
& -2 \sum_{i}(-1)^{i}\mathrm{dim}(H^{i}(T,\sO_{T})^{\Z_{m}}) \big)
\notag\end{split}\eeq
where $H^{i}(T,\cdot)^{\Z_{m}}$ is the subspace of $H^{i}(T,\cdot)$ fixed by $\Z_{m}$. The third equality holds since $\pi_{\ast}(\Omega_{T})^{\Z_{m}}=\Omega_{Y}^{[1]}$ and $\pi_{\ast}(\sO_{T})^{\Z_{m}}=\sO_{Y}$, where $\pi_{\ast}(\Omega_{T})^{\Z_{m}}$ and $\pi_{\ast}(\sO_{T})^{\Z_{m}}$ are the subsheaves of $\pi_{\ast}(\Omega_{T})$ and $\pi_{\ast}(\sO_{T})$, respectively, fixed by $\Z_{m}$.
Note that by Lemma \ref{Brauertr}
$$\mathrm{dim}(H^{i}(T,\Omega_{T})^{\Z_{m}})=\frac{1}{m}\sum_{g\in\Z_{m}}tr_{B}(g^{\ast}|_{H^{i}(\Omega_{T})}),\ \ \mathrm{dim}(H^{i}(T,\sO_{T})^{\Z_{m}})=\frac{1}{m}\sum_{g\in\Z_{m}}tr_{B}(g^{\ast}|_{H^{i}(\sO_{T})}).$$
Then by Lemma \ref{fixedpointformula}, we have
\beq\begin{split}
a(y,\Omega_{Y}^{[1]})=&\frac{1}{mN}\big(\sum_{i}(-1)^{i}\sum_{g\in\Z_{m}}tr_{B}(g^{\ast}|_{H^{i}(\Omega_{T})})-\chi(T,\Omega_{T})+2\ \chi(T,\sO_{T})\\
& -2 \sum_{i}(-1)^{i}\sum_{g\in\Z_{m}}tr_{B}(g^{\ast}|_{H^{i}(\sO_{T})}) \big)\\
=&\frac{1}{mN}\bigg(\sum_{i}(-1)^{i}\sum_{g\in\Z_{m}}tr_{B}(g^{\ast}|_{H^{i}(\Omega_{T})})-\big(\sum_{i}(-1)^{i}\sum_{g=\mathrm{id}}tr_{B}(g^{\ast}|_{H^{i}(\Omega_{T})})\big)\\
& +2\big(\sum_{i}(-1)^{i}\sum_{g=\mathrm{id}}tr_{B}(g^{\ast}|_{H^{i}(\sO_{T})})\big)-2 \sum_{i}(-1)^{i}\sum_{g\in\Z_{m}}tr_{B}(g^{\ast}|_{H^{i}(\sO_{T})}) \bigg)\\
=& \frac{1}{mN}\bigg(\sum_{i}(-1)^{i}\sum_{\mathrm{id}\neq g\in\Z_{m}}tr_{B}(g^{\ast}|_{H^{i}(\Omega_{T})})-2 \sum_{i}(-1)^{i}\sum_{\mathrm{id}\neq g\in\Z_{m}}tr_{B}(g^{\ast}|_{H^{i}(\sO_{T})}) \bigg)  \\
=&\frac{1}{mN}\big(N\sum_{\mathrm{id}\neq g\in\Z_{m}}\frac{tr_{B}(g|_{\Omega_{T,t}})}{(1-w(c_{g,1}))(1-w(c_{g,2}))}-2N\sum_{\mathrm{id}\neq g\in\Z_{m}}\frac{1}{(1-w(c_{g,1}))(1-w(c_{g,2}))}\big)\\
=&\frac{1}{m}\sum_{\mathrm{id}\neq g\in\Z_{m}}\frac{tr_{B}(g|_{\Omega_{T,t}})-2}{(1-w(c_{g,1}))(1-w(c_{g,2}))},
\notag\end{split}\eeq
where the fourth equality follows from (2) of Lemma \ref{projectivemodel}.
Then by \cite[Example 5.6 and Corollary 5.9]{langer2000chern},
$a(0,\Omega_{X}^{[1]})=a(y,\Omega_{Y}^{[1]})=\frac{1}{m}-1$, and hence the assertion holds.
\eproof
\etheorem

\bcorollary\label{aqoutientsing} Let $X$ and $f$ be as stated in Setting \ref{setting}.  We have 
$$c_{2}(f,\Omega_{\tilde{X}})\geq 2- \frac{1}{m}.$$
\bproof By Theorem \ref{localcontribution}, it suffices to show that $c_{2}(f,\Omega_{\tilde{X}})\geq 1- a(0,\Omega_{X}^{[1]})$.
Note that by Definition \ref{relative_Euler_characteristic},
$$c_{2}(f,\Omega_{\tilde{X}})=\chi(f,\Omega_{\tilde{X}})-2\chi(f,\sO_{\tilde{X}})+\frac{1}{2}c_{1}(f,\Omega_{\tilde{X}})\cdot(c_{1}(f,\Omega_{\tilde{X}})-K_{\tilde{X}})-a(0,\Omega_{X}^{[1]}).$$
Since $c_{1}(f,\Omega_{\tilde{X}})$ is supported on $E$ and $c_{1}(f,\Omega_{\tilde{X}})-K_{\tilde{X}}$ is $f$-trivial, we have
$$c_{1}(f,\Omega_{\tilde{X}})\cdot(c_{1}(f,\Omega_{\tilde{X}})-K_{\tilde{X}})=0.$$
Moreover, since $X$ has rational singularities, we have
$$\chi(f,\sO_{\tilde{X}})=h^{0}(f_{[\ast]}\sO_{\tilde{X}}/\sO_{X})+h^{0}(R^{1}f_{\ast}\sO_{\tilde{X}})=0.$$
Therefore, 
$$c_{2}(f,\Omega_{\tilde{X}})=\chi(f,\Omega_{\tilde{X}})-a(0,\Omega_{X}^{[1]}).$$
Since 
$$\chi(f,\Omega_{\tilde{X}})=h^{0}(f_{[\ast]}\Omega_{\tilde{X}}/f_{\ast}\Omega_{\tilde{X}})+h^{0}(R^{1}f_{\ast}\Omega_{\tilde{X}}),$$
we only need to show that $h^{0}(R^{1}f_{\ast}\Omega_{\tilde{X}})\geq 1$.
By the exact sequence
$$0\to\Omega_{\tilde{X}}(-E)\to\Omega_{\tilde{X}}\to\Omega_{\tilde{X}}|_{E}\to 0,$$
and $R^{2}f_{\ast}\Omega_{\tilde{X}}(-E)=0$ since $\mathrm{Exc}(f)$ is $1$-dimensional,
we know that the map 
$$R^{1}f_{\ast}\Omega_{\tilde{X}}\to R^{1}f_{\ast}(\Omega_{\tilde{X}}|_{E})\cong H^{1}(E,\Omega_{\tilde{X}}|_{E})$$
is surjective. Moreover, there exists a following exact sequence
$$\sO_{E}(-E)\to\Omega_{\tilde{X}}|_{E}\to\Omega_{E}\to 0.$$
It follows that the map 
$$H^{1}(E,\Omega_{\tilde{X}}|_{E})\to H^{1}(E,\Omega_{E})$$
is surjective. By Serre duality, $h^{1}(E,\Omega_{E})=h^{0}(E,\sO_{E})=1$.
Therefore  $h^{0}(R^{1}f_{\ast}\Omega_{\tilde{X}})\geq 1$, and hence the assertion holds.
\eproof
\ecorollary

\section{Proof of abundance when $\nu=2$ }

In this section, we prove the abundance conjecture in the case of $\nu(K_{X}+B)=2$. We first prove a special case (see Proposition \ref{specialcase}). Then we prove the general case by reducing it to the special case (see Theorem \ref{mainthm}).

\subsection{A key inequality}

\bsetting\label{setting_a_key_inequality} Let $(X,B)$ be an lc projective pair of dimension $3$ over an algebraically closed field $k$ of characteristic $p>3$. Assume that 

\noindent (1) $X$ is $\Q$-factorial, and $B$ is reduced, and $X\backslash B$ is terminal,

\noindent (2) there exists a divisor $D\in |m(K_{X}+B)|$ such that $D_{\mathrm{red}}=B$, where $m$ is a positive integer such that $L:=m(K_{X}+B)$ is Cartier,

\noindent (3) $K_{X}+B$ is nef and $\nu(K_{X}+B)=2$,

\noindent (4) if $C$ is a curve in $X$ with $C\cdot (K_{X}+B)>0$, then $(X,B)$ is plt at the generic  point of $C$.
\esetting

\blemma\label{2346} Let $S$ be a component of $B$ on which $L$ is not numerically trivial, and let $\lambda:S^{\lambda}\to S$ be the normalization. Let $\Theta$ be the $\Q$-divisor defined by adjunction
$(K_{X}+B)|_{S^{\lambda}}=K_{S^{\lambda}}+\Theta$. Then we have the following assertions hold.

\noindent (1) The Iitaka fibration $f$ associated to $K_{S^{\lambda}}+\Theta$ is a morphism to a smooth curve.

\noindent (2) Let $\Theta_{h}$ be the horizontal part of $\Theta$ with respect to $f$. We write $\Theta_{h}=\sum c_{q}\Gamma_{q}+\sum d_{l}\Delta_{l}$, where the $\Gamma_{q}$ map under $\lambda$ to the singular locus of $X$ and the $\Delta_{l}$ to the non-normal locus of $B$.
Then we have $f$ is a generically smooth map and one of the following holds.

\noindent (i)  The generic fibre of $\lambda$ is a smooth elliptic curve and $\Theta_{h}=0$,

\noindent (ii) the generic fibre of $\lambda$ is $\P^{1}$, $d_{l}=1$, and $c_{q}=1-\frac{1}{m_{q}}$ where $m_{q}$ are the indices of cyclic quotient singularities of $X$ at the generic point of $\Gamma_{q}$ (see Theorem \ref{lcadjunction}). Moreover, $(d_{1},\cdots;m_{1},\cdots)$ must be one of the following: $(1,1;)$, $(1;2)$, $(1;2,2)$, $(;2)$, $(;2,2)$, $(;2,2,2)$, $(;2,2,2,2)$, $(;2,3,6)$, $(;2,4)$, $(;2,4,4)$, $(;3)$, $(;3,3)$, $(;3,3,3)$.
\bproof For (1), note that 
$$0\leq (K_{S^{\lambda}}+\Theta)^{2} = S\cdot (K_{X}+B)^{2}=0
$$
where the last equality holds since by (2) and (3) in Setting \ref{setting_a_key_inequality}, $(K_{X}+B)^{2}$ is numerically trivial on $B$. Then $(S^{\lambda},\Theta)$ is a slc pair with $\nu(K_{S^{\lambda}}+\Theta)=1$ since $L$ is not numerically trivial on $S$.
It follows from Theorem \ref{slcabundance} that $K_{S^{\lambda}}+\Theta$ is semi-ample, and hence the Iitaka fibration $f$ associated to $K_{S^{\lambda}}+\Theta$ is a morphism to a smooth curve.
For (2), since $p>3$, by \cite[Proposition 2.2]{witaszek2021canonical}, $f$ is a generically smooth map. Let $F$ be the generic fibre of $f$ and consider the adjunction 
$$(K_{S^{\lambda}}+\Theta)|_{F}\sim_{\Q}K_{F}+\Theta_{h}|_{F}.$$
If $\Theta_{h}|_{F}=0$, then $\Theta_{h}=0$. We get (i). Otherwise, we have $\mathrm{deg}\ K_{F}<0$. Since $F$ is smooth, $F=\P^{1}$. It implies that $\mathrm{deg}\ \Theta_{h}|_{F}=2$. Now by (4) in Setting \ref{setting_a_key_inequality}, Proposition \ref{lcsurfacesing} and Theorem \ref{lcadjunction} we have (ii) holds.
\eproof
\elemma

The following proposition is the key inequality for proving Proposition \ref{specialcase}.

\bproposition\label{keyinequality}Let $\rho: V\to X$ be a desingularization. Then
$$\rho^{\ast}L \cdot\big(K^{2}_{V}+c_{2}(V)\big)\geq L\cdot\big(K_{X}^{2}+c_{2}(\Omega_{X}^{[1]})\big).$$
Moreover, if the equality holds, then we have $(K_{X}+B)|_{B^{\prime}}\sim_{\Q}K_{B^{\prime}}+\Delta_{B^{\prime}}$ such that $(K_{X}+B)|_{B^{\prime}}$ is numerically trivial on $\mathrm{Supp}\ \Delta_{B^{\prime}}$ where $B^{\prime}$ is the $S_{2}$-fication of $B$. 
\bproof By \cite[Corollary 2.2]{zhang2023frobenius} we can choose a sufficiently ample divisor $A$ such that $nL+A$ are very ample for all $n\gg 0$. We take a general section $H_{n}\in |nL+A|$. Note that $\rho^{\ast}L \cdot\big(K^{2}_{V}+c_{2}(V)\big)$ doesn't depend on the choice of $\rho$. Hence by Proposition \ref{gooddesingularization}, we can assume that $\rho: V\to X$ is a desingularization such that $\tilde{H}_{n}:=\rho^{-1}(H_{n})=\rho^{\ast}H_{n}$ is smooth, and $\rho|_{\tilde{H}_{n}}:\tilde{H}_{n}\to H_{n}$ is the minimal resolution of $H_{n}$.
Then by \cite[Theorem 0.4]{langer2025intersection} we have
\beq\begin{split}
\rho^{\ast}(nL+A) \cdot\big(K^{2}_{V}+c_{2}(V)\big)&=\tilde{H}_{n}\cdot\big(K^{2}_{V}+c_{2}(V)\big) \\
&= (K_{V}|_{\tilde{H}_{n}})^{2}+c_{2}(\Omega_{V}|_{\tilde{H}_{n}})
\notag\end{split}\eeq
and
\beq\begin{split}
(nL+A)\cdot\big(K_{X}^{2}+c_{2}(\Omega_{X}^{[1]})\big)&= H_{n}\cdot \big(K_{X}^{2}+c_{2}(\Omega_{X}^{[1]})\big)\\
&=(K_{X}|_{H_{n}})^{2}+c_{2}(\Omega^{[1]}_{X}|_{H_{n}}).
\notag\end{split}\eeq
We claim that
$$ \rho_{\ast}(K_{V}|_{\tilde{H}_{n}})^{2}+\rho_{\ast}c_{2}(\Omega_{V}|_{\tilde{H}_{n}})-(K_{X}|_{H_{n}})^{2}-c_{2}(\Omega^{[1]}_{X}|_{H_{n}})=\rho_{\ast}K_{\tilde{H}_{n}}^{2}+\rho_{\ast}c_{2}(\Omega_{\tilde{H}_{n}})-K_{H_{n}}^{2}-c_{2}(\Omega^{[1]}_{H_{n}})$$
as $0$-cycles.
Consider the exact sequence
\beq 0\to \sO_{\tilde{H}_{n}}(-\tilde{H}_{n})\to \Omega_{V}|_{\tilde{H}_{n}}\to \Omega_{\tilde{H}_{n}}\to 0. \eeq
Let $x$ be a point contained in the singular locus of $H_{n}$. Let $\rho_{x}$ be the map over $\mathrm{Spec}\ \sO_{H_{n},x}^{h}$ defined by (4.1).
By Proposition \ref{sumprinciple} and (1) of\cite[Proposition 3.2]{langer2025intersection}, we have
\beq\begin{split}
&\frac{1}{2}c_{1}(\rho_{x},\Omega_{V}|_{\tilde{H}_{n}})^{2}-c_{2}(\rho_{x},\Omega_{V}|_{\tilde{H}_{n}})\\
=&\frac{1}{2}c_{1}(\rho_{x},\Omega_{\tilde{H}_{n}})^{2}-c_{2}(\rho_{x},\Omega_{\tilde{H}_{n}})+\frac{1}{2}c_{1}(\rho_{x},\sO_{\tilde{H}_{n}}(-\tilde{H}_{n}))^{2}-c_{2}(\rho_{x},\sO_{\tilde{H}_{n}}(-\tilde{H}_{n}))\\
=&\frac{1}{2}c_{1}(\rho_{x},\Omega_{\tilde{H}_{n}})^{2}-c_{2}(\rho_{x},\Omega_{\tilde{H}_{n}})+\frac{1}{2}c_{1}(\rho_{x},\sO_{\tilde{H}_{n}}(-\tilde{H}_{n}))^{2}\\
=& \frac{1}{2}c_{1}(\rho_{x},\Omega_{\tilde{H}_{n}})^{2}-c_{2}(\rho_{x},\Omega_{\tilde{H}_{n}}),
\end{split}\eeq
where the last equality holds since 
$$\sO_{\tilde{H}_{n,x}}(\tilde{H}_{n,x})\sim \rho^{\ast}_{x}\sO_{H_{n}}(H_{n})$$
and by Definition \ref{relc1definition}
\beq c_{1}(\rho_{x},\sO_{\tilde{H}_{n}}(\tilde{H}_{n}))\sim_{\Q}\sO_{\tilde{H}_{n,x}}(\tilde{H}_{n,x})-\rho^{\ast}_{x}\sO_{H_{n}}(H_{n})\sim_{\Q}0.\eeq
By Definition \ref{relc1definition} and Definition \ref{c2definition}, we have
$$\sum_{x\in H_{n}\ \mathrm{singular}} c_{1}(\rho_{x},\Omega_{V}|_{\tilde{H}_{n}})=K_{V}|_{\tilde{H}_{n}}-\rho^{\ast} c_{1}(\Omega^{[1]}_{X}|_{H_{n}}),$$
$$\sum_{x\in H_{n}\ \mathrm{singular}} c_{2}(\rho_{x},\Omega_{V}|_{\tilde{H}_{n}})[x]=\rho_{\ast}c_{2}(\Omega_{V}|_{\tilde{H}_{n}})-c_{2}(\Omega^{[1]}_{X}|_{H_{n}}),$$
$$\sum_{x\in H_{n}\ \mathrm{singular}} c_{1}(\rho_{x},\Omega_{\tilde{H}_{n}})=K_{\tilde{H}_{n}}-\rho^{\ast} c_{1}(\Omega^{[1]}_{H_{n}})$$
and
$$\sum_{x\in H_{n}\ \mathrm{singular}}c_{2}(\rho_{x},\Omega_{\tilde{H}_{n}})[x]=\rho_{\ast}c_{2}(\Omega_{\tilde{H}_{n}})-c_{2}(\Omega^{[1]}_{H_{n}}).$$
They imply that 
\beq\begin{split}
&\rho_{\ast}(K_{V}|_{\tilde{H}_{n}})^{2}+\rho_{\ast}c_{2}(\Omega_{V}|_{\tilde{H}_{n}})-(K_{X}|_{H_{n}})^{2}-c_{2}(\Omega^{[1]}_{X}|_{H_{n}})\\
=&\rho_{\ast}\big(\rho^{\ast}(K_{X}|_{H_{n}})+\sum_{x\in H_{n}\ \mathrm{singular}} c_{1}(\rho_{x},\Omega_{V}|_{\tilde{H}_{n}})\big)^{2}-(K_{X}|_{H_{n}})^{2}+\sum_{x\in H_{n}\ \mathrm{singular}} c_{2}(\rho_{x},\Omega_{V}|_{\tilde{H}_{n}})[x]\\
=& \sum_{x\in H_{n}\ \mathrm{singular}} c_{1}(\rho_{x},\Omega_{V}|_{\tilde{H}_{n}})^{2}[x]+\sum_{x\in H_{n}\ \mathrm{singular}} c_{2}(\rho_{x},\Omega_{V}|_{\tilde{H}_{n}})[x]\\
=&\sum_{x\in H_{n}\ \mathrm{singular}}\frac{3}{2} c_{1}(\rho_{x},\Omega_{V}|_{\tilde{H}_{n}})^{2}[x]+\sum_{x\in H_{n}\ \mathrm{singular}} c_{2}(\rho_{x},\Omega_{V}|_{\tilde{H}_{n}})[x]-\sum_{x\in H_{n}\ \mathrm{singular}}\frac{1}{2} c_{1}(\rho_{x},\Omega_{V}|_{\tilde{H}_{n}})^{2}[x].
\notag\end{split}\eeq
Note that by (5.1), (5.3) and Definition \ref{relc1definition},
$$c_{1}(\rho_{x},\Omega_{V}|_{\tilde{H}_{n}})=c_{1}(\rho_{x},\Omega_{\tilde{H}_{n}})+c_{1}(\rho_{x},\sO_{\tilde{H}_{n}}(-\tilde{H}_{n}))=c_{1}(\rho_{x},\Omega_{\tilde{H}_{n}}).$$
Hence by (5.2) we have
\beq\begin{split}
&\rho_{\ast}(K_{V}|_{\tilde{H}_{n}})^{2}+\rho_{\ast}c_{2}(\Omega_{V}|_{\tilde{H}_{n}})-(K_{X}|_{H_{n}})^{2}-c_{2}(\Omega^{[1]}_{X}|_{H_{n}})\\
=&\sum_{x\in H_{n}\ \mathrm{singular}}\frac{3}{2} c_{1}(\rho_{x},\Omega_{\tilde{H}_{n}})^{2}[x]+\sum_{x\in H_{n}\ \mathrm{singular}} c_{2}(\rho_{x},\Omega_{\tilde{H}_{n}})[x]-\sum_{x\in H_{n}\ \mathrm{singular}}\frac{1}{2} c_{1}(\rho_{x},\Omega_{\tilde{H}_{n}})^{2}[x]\\
=&\sum_{x\in H_{n}\ \mathrm{singular}} c_{1}(\rho_{x},\Omega_{\tilde{H}_{n}})^{2}[x]+\sum_{x\in H_{n}\ \mathrm{singular}} c_{2}(\rho_{x},\Omega_{\tilde{H}_{n}})[x]\\
=&\rho_{\ast}\big(\rho^{\ast}K_{H_{n}}+\sum_{x\in H_{n}\ \mathrm{singular}} c_{1}(\rho_{x},\Omega_{\tilde{H}_{n}})\big)^{2}+\rho_{\ast}c_{2}(\Omega_{\tilde{H}_{n}})-K_{H_{n}}^{2}-c_{2}(\Omega^{[1]}_{H_{n}})\\
=&\rho_{\ast}K_{\tilde{H}_{n}}^{2}+\rho_{\ast}c_{2}(\Omega_{\tilde{H}_{n}})-K_{H_{n}}^{2}-c_{2}(\Omega^{[1]}_{H_{n}}).
\notag\end{split}\eeq
Hence, the claim holds.
Note that 
$$a_{n}:=\rho_{\ast}K_{\tilde{H}_{n}}^{2}+\rho_{\ast}c_{2}(\Omega_{\tilde{H}_{n}})-K_{H_{n}}^{2}-c_{2}(\Omega^{[1]}_{H_{n}})=\sum_{x\in H_{n}\ \mathrm{singular}} c_{1}(\rho_{x},\Omega_{\tilde{H}_{n}})^{2}[x]+\sum_{x\in H_{n}\ \mathrm{singular}} c_{2}(\rho_{x},\Omega_{\tilde{H}_{n}})[x]$$
is a $0$-cycle supported on the singular locus of $H_{n}$. Since $X\backslash B$ is terminal by (1) in Setting \ref{setting_a_key_inequality}, $a_{n}$ is supported on $H_{n}\cap C$ where $C:=C_{1}+C_{2}$, $C_{1}\subseteq B$ is the union of curves contained in the singular locus of $X$ on which $L$ are numerically trivial and $C_{2}\subseteq B$ is the union of curves contained in the singular locus of $X$ on which $L$ are positive.

We claim that the degrees of $a_{n}$ on $C_{1}$ are bounded from below.
Note that for any irreducible component $\Gamma$ of $C_{1}$, 
$$\Gamma\cdot H_{n}=\Gamma\cdot (nL+A)=\Gamma\cdot A$$
are independent of $n$.
Let $z_{1,n}$ be a point in $H_{n}\cap\Gamma$. It suffices to show that the degrees of $a_{n}$ at $z_{1,n}$ are bounded from below.
By Theorem \ref{bertiniforklt}, the dual graphs of $(H_{n},z_{1,n})$ are independent of $n$.  Hence, $c_{1}(\rho_{z_{1,n}},\Omega_{\tilde{H}_{n}})^{2}$
are bounded from below since $\rho|_{\tilde{H}_{n}}$ are minimal resolutions and $c_{1}(\rho_{z_{1,n}},\Omega_{\tilde{H}_{n}})^{2}$ are determined by the dual graphs of $(H_{n},z_{1,n})$.
By \cite[Proposition 3.2]{langer2025intersection}, 
$$c_{1}(\rho_{z_{1,n}},\Omega_{\tilde{H}_{n}})^{2}+c_{2}(\rho_{z_{1,n}},\Omega_{\tilde{H}_{n}})\geq \frac{5}{4}c_{1}(\rho_{z_{1,n}},\Omega_{\tilde{H}_{n}})^{2}.$$
Hence the claim follows. Therefore, to prove that
$\rho^{\ast}L \cdot\big(K^{2}_{V}+c_{2}(V)\big)- L\cdot\big(K_{X}^{2}+c_{2}(\Omega_{X}^{[1]})\big)\geq 0,$
it suffices to show that
$(nL+A) \cdot\big(\rho_{\ast}(K^{2}_{V}+c_{2}(V))-(K_{X}^{2}+c_{2}(\Omega_{X}^{[1]}))\big)=a_{n}$
are of non-negative degrees on $C_{2}$ for all $n\gg 0$. Fix an $n\gg 0$.  
Let $z_{2}\in C_{2}$ be a singular point of $H_{n}$. 
Note that $z_{2}$ is in an irreducible component $S$ of $B$, on which $L$ is not numerically trivial. Then we can apply Lemma \ref{2346} and Corollary \ref{bertiniforcyclicquotient} to conclude that
$(H_{n},z_{2})$ is a cyclic quotient singularity whose index $m$ is one of $2,3,4,6$. Then it follows from Corollary \ref{aqoutientsing} that, if $z_{2}$ is a Du Val singularity, 
the degree of $a_{n}$ at $z_{2}$, which is equal to $c_{2}(\rho_{z_{2}},\Omega_{\tilde{H}_{n}})$, is at least $\frac{3}{2},\frac{5}{3},\frac{7}{4},\frac{11}{6}$ when $m$ is $2,3,4,6$ respectively. Otherwise, 
by \cite[Definition 3.33]{kollar2013singularities},
$m$ is one of $3,4,6$, and $c_{1}(\rho_{z_{2}},\Omega_{\tilde{H}_{n}})^{2}$ is $-\frac{1}{3}, -1, -\frac{8}{3}$ when $m$ is $3,4,6$ respectively.
Hence by Corollary \ref{aqoutientsing}, 
the degree of $a_{n}$ at $z_{2}$, is at least $\frac{4}{3},\frac{3}{4},-\frac{5}{6}$ when $m$ is $3,4,6$ respectively. By Lemma \ref{2346} and Corollary \ref{bertiniforcyclicquotient}
, an index $6$ point is always accompanied by an index $2$ point and an index $3$ point. Moreover, different index $6$ points give different index $2$ points and index $3$ points. Hence we always have $a_{n}$ is of positive degree on $C_{2}$.
Therefore, the inequality holds. Moreover, if the equality holds, then there is no such singular point $z_{2}$. It is to say that (ii) holds in (2) of Lemma \ref{2346}. Hence, the second assertion holds by Lemma \ref{slcsurface}. 
\eproof
\eproposition

\subsection{A special case}
In this subsection, we prove a special case of the abundance conjecture in the case of $\nu(K_{X}+B)=2$ (see Proposition \ref{specialcase}).
Before proving Proposition \ref{specialcase}, we give some results which help us to calculate some cohomology.

\blemma\label{Bcohomology} Let $B$ be a reduced equidimensional quasi-projective scheme of dimension $2$ over a field and $g:B^{\prime}\to B$ be its $S_{2}$-fication. Assume that $L$ is a Cartier divisor on $B$. Then 
$$h^{i}(B^{\prime},\sO_{B^{\prime}}(ng^{\ast}L))-h^{i}(B,\sO_{B}(nL))$$
 are bounded for all $i$ and $n\geq 0$.
\bproof Since $g$ is a finite morphism by Proposition \ref{S2fication}, we have
$$h^{i}(B^{\prime}, \sO_{B^{\prime}}(ng^{\ast}L))=h^{i}(B,g_{\ast}\sO_{B^{\prime}}\otimes \sO_{B}(nL)).$$
Note that there exists an exact sequence
$$0\to\sO_{B}\to g_{\ast}\sO_{B^{\prime}}\to \sQ\to 0$$
such that $\sQ$ is a coherent sheaf supported on some points by Proposition \ref{S2fication}.
Tensoring it with $\sO_{B}(nL)$, we get
$$0\to\sO_{B}(nL)\to g_{\ast}\sO_{B^{\prime}}\otimes\sO_{B}(nL)\to \sQ\otimes\sO_{B}(nL)\to 0.$$
They give the exact sequences
$$H^{i-1}(B,\sQ\otimes\sO_{B}(nL))\to H^{i}(B,\sO_{B}(nL))\to H^{i}(B,g_{\ast}\sO_{B^{\prime}}\otimes\sO_{B}(nL))\to H^{i}(B,\sQ\otimes\sO_{B}(nL)).$$
Since both $h^{i-1}(B,\sQ\otimes\sO_{B}(nL))$ and $h^{i}(B,\sQ\otimes\sO_{B}(nL))$ are constants for $n\geq 0$, the assertion holds.
\eproof
\elemma

\blemma \label{vanishing} Let $X$ be a klt projective threefold over an algebraically closed field of characteristic $p>3$. Assume that $M$ is a nef $\Q$-Cartier Weil divisor on $X$ such that $\nu(M)\geq 2$. Then we have 
$h^{2}(X,\sO_{X}(K_{X}+nM))$
is bounded for $n\geq 0$. 
\bproof
We take a general hyperplane section $A$ on $X$ and
consider the exact sequences
$$0\to\sO_{X}(K_{X}+nM+lA)\to \sO_{X}(K_{X}+nM+(l+1)A)\to \sQ_{n,l+1}\to 0, $$
where $l\geq 0$ is an integer. Since both  $\sO_{X}(K_{X}+nM+lA)$ and $\sO_{X}(K_{X}+nM+(l+1)A)$ are $S_{2}$, by \cite[Lemma 2.60]{kollar2013singularities} $\sQ_{n,l+1}$ are $S_{1}$ and supported on $A$. Hence they are torsion-free. Note that the sheaves $\sO_{A}\Big(\big(K_{X}+nM+(l+1)A\big)|_{A}\Big)$ are reflexive. They give exact sequences
$$0\to\sQ_{n,l+1}\to \sO_{A}\Big(\big(K_{X}+nM+(l+1)A\big)|_{A}\Big)\to \sG_{n,l+1}\to 0$$
where $\sG_{n,l+1}$ are supported on some points. Moreover, since $M$ is $\Q$-Cartier, $\sG_{n,l+1}$ are of bounded ranks at every point.
Hence there are exact sequences
$$H^{0}(A,\sG_{n,l+1})\to H^{1}(A,\sQ_{n,l+1})\to H^{1}\bigg(A,\sO_{A}\Big(\big(K_{X}+nM+(l+1)A\big)|_{A}\Big)\bigg).$$
Since $A$ is a klt surface by Theorem \ref{bertiniforklt}, for $n\gg 0$ and any fixed $l$ 
$$H^{1}\bigg(A,\sO_{A}\Big(\big(K_{X}+nM+(l+1)A\big)|_{A}\Big)\bigg)=0$$
by \cite[Theorem 2.11]{tanaka2012x}.
Therefore, both $h^{0}(A,\sG_{n,l+1})$ and $h^{1}\bigg(A,\sO_{A}\Big(\big(K_{X}+nM+(l+1)A\big)|_{A}\Big)\bigg)$ are bounded for $n\geq 0$ and any fixed $l$.
It implies that there exists a constant $C_{l+1}$ such that 
$$h^{1}(X,\sQ_{n,l+1})\leq C_{l+1}$$ for $n\geq 0$ and any $l$.
By \cite[Theorem 1.5]{keeler2003ample}, there exists an integer $l_{0}> 0$ such that for $l\geq l_{0}$ and all $n\geq 0$ we have $h^{2}\Big(X,\sO_{X}\big(K_{X}+nM+(l+1)A\big)\Big)=0$.
It follows that
$$h^{2}(X,\sO_{X}(K_{X}+nM))\leq h^{1}(X,\sQ_{n,1})+h^{2}(X,\sO_{X}(K_{X}+nM+A))\leq \cdots\leq \sum_{i=1}^{l_{0}+1}h^{1}(X,\sQ_{n,i})\leq\sum_{i=1}^{l_{0}+1}C_{i}.$$

\eproof
\elemma

\bproposition\label{specialcase} Let $(X,B)$ be an lc projective pair of dimension $3$ over an algebraically closed field $k$ of characteristic $p>3$. Assume that 

\noindent (1) $X$ is $\Q$-factorial, $B$ is reduced, and $X\backslash B$ is terminal,

\noindent (2) there exists a divisor $D\in |m(K_{X}+B)|$ such that $D_{\mathrm{red}}=B$ where $m$ is a positive integer such that $L:=m(K_{X}+B)$ is Cartier,

\noindent (3) $K_{X}+(1-\varepsilon)B$ is nef for any sufficiently small $\varepsilon \geq 0$,

\noindent (4)  $\nu(K_{X}+B)=2$,

\noindent (5) if $C$ is a curve in $X$ with $C\cdot (K_{X}+B)>0$, then $(X,B)$ is plt at the generic  point of $C$.

\noindent Then $K_{X}+B$ is semi-ample.
\bproof By Theorem \ref{known} it suffices to show $\kappa(K_{X}+B)> 0$. Since $K_{X}+(1-\varepsilon)B$ is nef for any sufficiently small $\varepsilon\geq 0$, by  Theorem \ref{nonvan} there exists an effective $\Q$-divisor $\Delta_{\varepsilon}\sim_{\Q}K_{X}+(1-\varepsilon)B$ for any sufficiently small rational $\varepsilon\geq 0$. Moreover, by (2) we can assume that $B\subseteq \mathrm{Supp}\ \Delta_{\varepsilon}$ and $\Delta_{\varepsilon}$ have the same support for all sufficiently small rational $\varepsilon\geq 0$. Therefore by \cite[Lemma 11.3.3]{kollar1992flips} we have
$$\nu(K_{X}+(1-\varepsilon)B)=\nu(K_{X}+B)=2$$
for any sufficiently small rational $\varepsilon>0$.
It implies 
$$(K_{X}+(1-\varepsilon)B)^{3}=(K_{X}+B)^{3}-3\varepsilon(K_{X}+B)^{2}\cdot B+3\varepsilon^{2}(K_{X}+B)\cdot B^{2}-\varepsilon^{3}B^{3}=0.$$
Hence we have 
\beq (K_{X}+B)^{3}=K_{X}\cdot (K_{X}+B)^{2}=K_{X}^{2}\cdot (K_{X}+B)=B^{3}=0.\eeq

Now let $\rho:V\to X$ be a desingularization of $X$. By Riemann-Roch theorem,
$$\chi(V,\sO_{V}(n\rho^{\ast}L))=\frac{n^{3}}{6}(\rho^{\ast}L)^{3}-\frac{n^{2}}{4}K_{V}\cdot(\rho^{\ast}L)^{2}+\frac{n}{12}(c_{2}(V)+K_{V}^{2})\cdot \rho^{\ast}L+\chi(\sO_{V}).$$ 
Note that by the projection formula and (5.4),
$$(\rho^{\ast}L)^{3}=K_{V}\cdot(\rho^{\ast}L)^{2}=0.$$
Hence we have
$$\chi(V,\sO_{V}(n\rho^{\ast}L))=\frac{n}{12}(K_{V}^{2}+c_{2}(V))\cdot \rho^{\ast}L+\chi(\sO_{V}).$$

If $T_{V}$ is not strongly $(\rho^{\ast}L,\rho^{\ast}L)$-semistable, then by Theorem \ref{constructfoliation} there is a $1$-foliation $\sF\subseteq T_{V}$ such that $\mu(\sF)>0=\mu(T_{V})$ w.r.t. $(\rho^{\ast}L,\rho^{\ast}L)$. By \cite[Corollary 2.2]{zhang2023frobenius} we can take a sufficiently ample line bundle $H$ on $V$ such that $a\rho^{\ast}L+H$ are very ample for all integers $a\geq 0$. Then we can use two general divisors in $|a\rho^{\ast}L+H|$ to cut out a smooth curve $C_{a}$ such that $\sF$ is smooth along $C_{a}$.
Note that
\beq\begin{split}
c_{1}(\sF)\cdot C_{a}&=c_{1}(\sF)\cdot(a\rho^{\ast}L+H)^{2}\\
&=a^{2}c_{1}(\sF)\cdot(\rho^{\ast}L)^{2} +2a c_{1}(\sF)\cdot\rho^{\ast}L\cdot H+c_{1}(\sF)\cdot H^{2},    \notag
\end{split}\eeq
\beq\begin{split}
K_{V}\cdot C_{a}&=K_{V}\cdot(a\rho^{\ast}L+H)^{2}\\
&=2a K_{V}\cdot\rho^{\ast}L\cdot H+K_{V}\cdot H^{2}.    \notag
\end{split}\eeq
and 
\beq\begin{split}
\rho^{\ast}L\cdot C_{a}&=\rho^{\ast}L\cdot(a\rho^{\ast}L+H)^{2}\\
&=2a (\rho^{\ast}L)^{2}\cdot H+\rho^{\ast}L\cdot H^{2}.    \notag
\end{split}\eeq
Since $\mu(\sF)=\frac{c_{1}(\sF)\cdot(\rho^{\ast}L)^{2}}{\mathrm{rk}\sF}>0$, 
$c_{1}(\sF)\cdot C_{a}$ grows with $a^2$, while $K_{V}\cdot C_{a}$ and $\rho^{\ast}L\cdot C_{a}$ grow at most linearly in $a$. In particular, for $a\gg 0$
$$c_{1}(\sF)\cdot C_{a} > \frac{K_{V}\cdot C_{a}}{p-1}.$$
By Theorem \ref{bendbreakfoliation}, we know for general $y\in C_{a}$ there exists a rational curve $B_{y}$ such that $B_{y}$ passes through $y$ and
$$\rho^{\ast}L\cdot B_{y}\leq \frac{6p\rho^{\ast}L\cdot C_{a}}{(p-1)c_{1}(\sF)\cdot C_{a}-K_{V}\cdot C_{a}}.$$
Note that $\rho^{\ast}L\cdot B_{y}$ is a non-negative integer, while the denominator of the RHS above grows with $a^2$, and the numerator grows at most linearly in $a$. 
Hence $\rho^{\ast}L\cdot B_{y}=0$ for $a\gg 0$. Since $C_{a}$ and $y$ are general, 
we obtain a family of curves $B_{y}$ covering $V$ such that $\rho^{\ast}L\cdot B_{y}=0$.
It follows that
$$n(K_{X}+B)=n(L)=n(\rho^{\ast}L)\leq 2.$$
By Theorem \ref{known}, we have $K_{X}+B$ is semi-ample.

From now on, we can assume that $T_{V}$ is strongly  $(\rho^{\ast}L,\rho^{\ast}L)$-semistable.
Then we have $\Omega^{[1]}_{X}$ is strongly  $(L,L)$-semistable since $\rho$ is birational. Hence by Theorem \ref{Bogomolovineq}, 
$$\Delta(\Omega^{[1]}_{X})\cdot L=\big(6c_{2}(\Omega^{[1]}_{X})-2c_{1}(\Omega^{[1]}_{X})^{2}\big)\cdot L\geq 0.$$
Since $K_{X}^{2}\cdot L=K_{X}^{2}\cdot (K_{X}+B)=0$ by (5.4), we have
\beq c_{2}(\Omega^{[1]}_{X})\cdot L\geq 0.\eeq
Note that by (1)-(5), the assumptions in Setting \ref{setting_a_key_inequality} are satisfied.
Hence we can apply Proposition \ref{keyinequality} and (5.4) to conclude that
\beq\begin{split}
\chi(V,\sO_{V}(n\rho^{\ast}L))&=\frac{n}{12}(K_{V}^{2}+c_{2}(V))\cdot \rho^{\ast}L+\chi(\sO_{V})\\
&\geq \frac{n}{12}(K_{X}^{2}+c_{2}(\Omega^{[1]}_{X}))\cdot L+\chi(\sO_{V})\\
&= \frac{n}{12}c_{2}(\Omega^{[1]}_{X})\cdot L+\chi(\sO_{V}).
\end{split}\eeq
We claim that there exists a constant $t_{1}$ such that 
\beq\chi(X,\sO_{X}(nL))\geq \chi(V,\sO_{V}(n\rho^{\ast}L))+t_{1}.\eeq
Note that $R^{i}\rho_{\ast}\sO_{V}$ are supported on the closed subsets of $X$ which are of dimension at most $2-i$ for $i=1,2$.
By Leray spectral sequence, we have
\beq\begin{split}
&\chi(V,\sO_{V}(n\rho^{\ast}L))\\
=&h^{0}(V,\sO_{V}(n\rho^{\ast}L))-h^{1}(V,\sO_{V}(n\rho^{\ast}L))+h^{2}(V,\sO_{V}(n\rho^{\ast}L))-h^{3}(V,\sO_{V}(n\rho^{\ast}L))\\
=&h^{0}(X,\sO_{X}(nL))-\big(h^{1}(X,\sO_{X}(nL))+h^{0}(X,R^{1}\rho_{\ast}\sO_{V}\otimes\sO_{X}(nL))\big)+\big(h^{2}(X,\sO_{X}(nL))+\\
&\ h^{1}(X,R^{1}\rho_{\ast}\sO_{V}\otimes\sO_{X}(nL))+h^{0}(X,R^{2}\rho_{\ast}\sO_{V}\otimes\sO_{X}(nL))\big)-h^{3}(X,\sO_{X}(nL))\\
=&\chi(X,\sO_{X}(nL))-h^{0}(X,R^{1}\rho_{\ast}\sO_{V}\otimes\sO_{X}(nL))+h^{1}(X,R^{1}\rho_{\ast}\sO_{V}\otimes\sO_{X}(nL))+h^{0}(X,R^{2}\rho_{\ast}\sO_{V}\otimes\sO_{X}(nL)).
\notag\end{split}\eeq
It follows that 
$$\chi(X,\sO_{X}(nL))+h^{1}(X,R^{1}\rho_{\ast}\sO_{V}\otimes\sO_{X}(nL))+h^{0}(X,R^{2}\rho_{\ast}\sO_{V}\otimes\sO_{X}(nL))\geq \chi(V,n\rho^{\ast}L)+h^{0}(X,R^{1}\rho_{\ast}\sO_{V}\otimes\sO_{X}(nL)).$$
Note that $h^{0}(X,R^{2}\rho_{\ast}\sO_{V}\otimes\sO_{X}(nL))$ is bounded for $n\geq 0$ since $R^{2}\rho_{\ast}\sO_{V}$ are supported on some points. Since $R^{1}\rho_{\ast}\sO_{V}$ is a coherent sheaf supported on a closed subset $Z_{1}$ of $X$ which is of dimension at most $1$, $h^{1}(X,R^{1}\rho_{\ast}\sO_{V}\otimes\sO_{X}(nL))$ is bounded for $n\geq 0$ by \cite[Theorem 1.4.40]{lazarsfeld2017positivity}. Hence (5.7) holds.

Now let's consider the exact sequences
\beq 0\to \sO_{X}(nL-B)\to\sO_{X}(nL)\to\sO_{B}(nL|_{B})\to 0.\eeq
We claim that $h^{2}(X,nL)$ is bounded for $n\geq 0$.
We write $g:B^{\prime}\to B$ for the $S_{2}$-fication of $B$. 
By Lemma \ref{slcsurface}, we have $(K_{X}+B)|_{B^{\prime}}\sim_{\Q}K_{B^{\prime}}+\Delta_{B^{\prime}}$
and $(B^{\prime},\Delta_{B^{\prime}})$ is a slc pair, where $|_{B'}$ denotes the pullback to $B'$.
Note that by (5.4) 
$$(K_{B^{\prime}}+\Delta_{B^{\prime}})^{2}=B\cdot (K_{X}+B)^{2}=0$$
and by (2) and (4)
$$(K_{B^{\prime}}+\Delta_{B^{\prime}})\cdot A|_{B^{\prime}}=B\cdot (K_{X}+B)\cdot A>0$$
where $A$ is an ample divisor on $X$.
Hence $\nu(K_{B^{\prime}}+\Delta_{B^{\prime}})=1$. 
By Theorem \ref{slcabundance}, 
$$\kappa(L|_{B'})=\kappa(K_{B^{\prime}}+\Delta_{B^{\prime}})=\nu(K_{B^{\prime}}+\Delta_{B^{\prime}})=1.$$
Hence 
$$h^{0}\Big(B^{\prime},\sO_{B^{\prime}}\big(nm(K_{B^{\prime}}+\Delta_{B^{\prime}})\big)\Big)=h^{0}(B^{\prime},\sO_{B^{\prime}}(n L|_{B'}))$$
grows with $n$.
By Serre duality,
\beq\begin{split}
h^{2}(B^{\prime},\sO_{B^{\prime}}(n L|_{B'}))&=h^{0}(B^{\prime},\sO_{B^{\prime}}(K_{B^{\prime}}-n L|_{B'}))\\
&=h^{0}\Big(B^{\prime},\sO_{B^{\prime}}\big(K_{B^{\prime}}-nm(K_{X}+B)|_{B^{\prime}}\big)\Big)\\
&= h^{0}\Big(B^{\prime},\sO_{B^{\prime}}\big(K_{B^{\prime}}-nm(K_{B^{\prime}}+\Delta_{B^{\prime}})\big)\Big)\\
&=0   \notag
\end{split}\eeq
for $n\gg 0$. It follows that 
$$h^{2}(B,\sO_{B}(nL|_{B}))$$
 is bounded for $n\geq 0$ by Lemma \ref{Bcohomology}. Moreover, since 
$$nL-B\sim_{\Q}K_{X}+(nm-1)(K_{X}+B)$$
 and $K_{X}+B$ is a nef Weil divisor such that $\nu(K_{X}+B)=2$, by Lemma \ref{vanishing} we have
$$h^{2}(X,\sO_{X}(nL-B))$$
 is bounded for $n\geq 0$. Therefore, by (5.8), $h^{2}(X,\sO_{X}(nL))$ is bounded for $n\geq 0$.

If the inequality in (5.6) is strict, then by (5.5)
$$(K_{V}^{2}+c_{2}(V))\cdot \rho^{\ast}L > c_{2}(\Omega^{[1]}_{X})\cdot L\geq 0.$$
Hence $\chi(V,\sO_{V}(n\rho^{\ast} L))$ grows with $n$.
By (5.7), $\chi(X,\sO_{X}(nL))$, and hence $h^{0}(X,\sO_{X}(nL))$, grow with $n$. It implies $\kappa(K_{X}+B)>0$. The assertion holds.
Otherwise, we have $(K_{X}+B)|_{B^{\prime}}$ is numerically trivial on $\mathrm{Supp}\ \Delta_{B^{\prime}}$ by Proposition \ref{keyinequality}.
Since 
$$h^{2}(X,\sO_{X}(nL))$$ 
is bounded for $n\geq 0$, by (5.5), (5.6) and (5.7) 
there exists a constant $t_{2}$ such that
\beq h^{0}(X,\sO_{X}(nL))\geq h^{1}(X,\sO_{X}(nL))+t_{2}.\eeq
Note that by Theorem \ref{RRslc} and (5.4)
\beq \begin{split}
\chi(B^{\prime},\sO_{B^{\prime}}(n L|_{B'}))&=\frac{1}{2}(nL|_{B^{\prime}})\cdot(nL|_{B^{\prime}}-K_{B^{\prime}})+\chi(B^{\prime},\sO_{B^{\prime}})\\
&= \frac{1}{2}(nL|_{B^{\prime}})\cdot(-K_{B^{\prime}}-\Delta_{B^{\prime}}+\Delta_{B^{\prime}})+\chi(B^{\prime},\sO_{B^{\prime}})\\
&= \frac{1}{2}(nL|_{B^{\prime}})\cdot\Delta_{B^{\prime}}+\chi(B^{\prime},\sO_{B^{\prime}})\\
&= \chi(B^{\prime},\sO_{B^{\prime}}). \notag
\end{split}\eeq
Since $h^{0}(B^{\prime},\sO_{B^{\prime}}(n L|_{B'}))$ grows with $n$ and $h^{2}(B^{\prime},\sO_{B^{\prime}}(n L|_{B'}))$ is bounded for $n\geq 0$, we have $h^{1}(B^{\prime},\sO_{B^{\prime}}(n L|_{B'}))$ grows with $n$. By Lemma \ref{Bcohomology}, 
$$h^{1}(B,\sO_{B}(n L|_{B}))$$ 
grows with $n$.
Since 
$$h^{2}(X,\sO_{X}(nL-B))$$ 
is bounded for $n\geq 0$, by (5.8), 
$$h^{1}(X,\sO_{X}(nL))$$ 
grows with $n$.
Hence, by (5.9), $h^{0}(X,\sO_{X}(nL))$ grows with $n$.
In conclusion, we have 
$$\kappa(K_{X}+B)=\kappa(X,L)> 0,$$ 
and hence $K_{X}+B$ is semi-ample.
\eproof
\eproposition

\subsection{The general case}

We recall a standard lemma from \cite[Lemma 14.2]{kollar1992flips}.  

\blemma\label{standardmodification}  Let $(X,B)$ be a $\Q$-factorial dlt projective pair of dimension $3$ over an algebraically closed field $k$ of characteristic $p>3$. Assume that $X$ is terminal, $K_{X}+B$ is nef and there exists an effective $\Q$-divisor $D$ such that
$D\sim_{\Q} K_{X}+B$ and $\mathrm{Supp}\ B\subseteq\mathrm{Supp}\ D$. Then there exists a $\Q$-factorial lc pair $(Y,B_{Y})$ such that

\noindent (1) $(Y\backslash\mathrm{Supp}\ B_{Y})\cong (X\backslash\mathrm{Supp}\  D)$ is terminal, and $B_{Y}$ is reduced,

\noindent (2) $K_{Y}+B_{Y}\sim_{\Q} D_{Y}$ for an effective $\Q$-divisor $D_{Y}$ with $\mathrm{Supp}\ D_{Y}=  B_{Y},$

\noindent (3) $K_{Y}+(1-\varepsilon)B_{Y}$ is nef for any sufficiently small $\varepsilon>0$,

\noindent (4) $\kappa(K_{X}+B)= \kappa(K_{Y}+B_{Y})$ and $\nu(K_{X}+B)= \nu(K_{Y}+B_{Y})$,

\noindent (5) if $C$ is a curve in $X$ with $C\cdot (K_{X}+B)>0$, then $(X,B)$ is plt at the generic  point of $C$.
\elemma
\bproof Let $g: W\to X$ be a log smooth resolution of $(X,\mathrm{Supp}\ D)$ with the reduced exceptional divisor $E$. Set $B_{W}:= \mathrm{Supp}(g^{-1}_{\ast}D)+ E$.
Since $(X,B)$ is dlt, we can write $K_{W}+B_{W}\sim_{\Q} D_{W}$, where
$D_{W}= g^{\ast}D+(K_{W}+B_{W}-g^{\ast}(K_{X}+B))$
is an effective $\Q$-divisor with 
$$\mathrm{Supp}\ D_{W}= \mathrm{Supp}(g^{-1}_{\ast}D)+E= B_{W}.$$
 By replacing $W$ by the output of a $(K_{W}+B_{W})$-MMP over $X$ by Theorem \ref{lcmmp} and using again that $(X,B)$ is dlt, we can assume that $W \backslash \mathrm{Supp}\ D_{W}\cong X \backslash \mathrm{Supp}\ D$ and $\mathrm{Supp}\ D_{W}= \mathrm{Supp}\  g^{\ast}D$.
 
Now, since $K_{W}+B_{W}$ is pseudoeffective, by Theorem \ref{lcmmp} we can run a $(K_{W} +B_{W})$-MMP which terminates with a minimal model $h: W \dashrightarrow Y$. In particular $K_{Y}+B_{Y}$ is nef for $B_{Y}:= h_{\ast}B_{W}$. Furthermore, $K_{Y}+B_{Y}\sim_{\Q} D_{Y}$ and $\mathrm{Supp}\ D_{Y} = B_{Y}$, where $D_{Y} := h_{\ast}D_{W}$. Hence, (2) holds. Since $h$ and $g$ are
isomorphisms outside of $D_{W}$ , we get (1). Note that both $K_{X}+B$ and $K_{Y}+B_{Y}$ are effective and nef $\Q$-divisors.  Applying \cite[Lemma 11.3.3]{kollar1992flips} to pullbacks of  $K_{X}+B$ and $K_{Y}+B_{Y}$  on any common resolution of $X$ and $Y$, we get
$$\kappa(K_{X}+B)= \kappa(K_{Y}+B_{Y}),\nu(K_{X}+B)= \nu(K_{Y}+B_{Y}).$$
Hence (4) holds. To get (3), note that since $K_{Y}+B_{Y}\sim_{\Q}D_{Y}$ and $\mathrm{Supp}\ D_{Y}= B_{Y}$, we have $K_{Y}+(1-\varepsilon)B_{Y}$ is effective for any sufficiently small $\varepsilon\geq 0$. Hence by Theorem \ref{Waldrontermination}, we can run a $K_{Y}$-MMP which is $(K_{Y}+B_{Y})$-trivial and it terminates. Replacing $(Y,B_{Y})$ by the output, we get (3).
Finally, since cyclic quotient surface singularities are klt (see \cite[Theorem A.3]{sato2025general} for example), (5) follows from Lemma \ref{2346} and the proof of \cite[Lemma 14.2]{kollar1992flips}.
\eproof

Now we can prove our main result.

\btheorem\label{mainthm} Let $(X,B)$ be an lc projective pair of dimension $3$ over a perfect field $k$ of characteristic $p>3$. Assume that $K_{X}+B$ is nef and $\nu(K_{X}+B)=2$. Then $K_{X}+B$ is semi-ample.
\bproof By \cite[Remark 2.7]{gongyo2019rational}, we may assume that $k$ is algebraically closed. By Theorem \ref{dltmodification}, replacing $(X,B)$ by a dlt model, we can assume that $X$ is terminal, and $(X,B)$ is $\Q$-factorial and dlt. By Definition \ref{MMPKtrivial} we run a $K_{X}$-MMP which is $(K_{X}+B)$-trivial. Since $X$ is terminal, it terminates by Lemma \ref{terminaltermination}. If it terminates with a Mori fibre space, then we have $n(K_{X}+B)\leq 2$. By Theorem \ref{known}, $K_{X}+B$ is semi-ample. Hence by Lemma \ref{Output} we can assume that the $K_{X}$-MMP which is $(K_{X}+B)$-trivial terminates with an lc pair $(X^{\prime},B^{\prime})$ such that $K_{X^{\prime}}+(1-\varepsilon)B^{\prime}$ is nef  for any sufficiently small $\varepsilon>0$. By Theorem \ref{nonvan}, we have
$$\kappa(K_{X}+(1-\varepsilon)B)=\kappa(K_{X^{\prime}}+(1-\varepsilon)B^{\prime})\geq 0$$ 
 for any sufficiently small rational $\varepsilon>0$. Hence there exists an effective $\Q$-divisor $\Delta\sim_{\Q}K_{X}+B$ such that $\mathrm{Supp}\ B\subseteq\mathrm{Supp}\ \Delta$. 
 By Lemma \ref{standardmodification}, we can reduce the assertion to the case when $(X,B)$ is an lc projective pair of dimension $3$ and the following assertions hold, 

\noindent (1) $X$ is $\Q$-factorial, $B$ is reduced, and $X\backslash B$ is terminal,

\noindent (2) $K_{X}+B\sim_{\Q} D$ for an effective $\Q$-divisor $D$ with $\mathrm{Supp}\ D=  B,$

\noindent (3) $K_{X}+(1-\varepsilon)B$ is nef for any sufficiently small $\varepsilon>0$,

\noindent (4)  $\nu(K_{X}+B)=2$,

\noindent (5) if $C$ is a curve in $X$ with $C\cdot (K_{X}+B)>0$, then $(X,B)$ is plt at the generic  point of $C$.

\noindent Then by Proposition \ref{specialcase}, we have $K_{X}+B$ is semi-ample.
\eproof
\etheorem

\bibliographystyle{plain}
 \bibliography{Reference_CM9004_final}

@article{arvidsson2022kawamata,
  title={On the Kawamata--Viehweg vanishing theorem for log del Pezzo surfaces in positive characteristic},
  author={Arvidsson, Emelie and Bernasconi, Fabio and Lacini, Justin},
  journal={Compositio Mathematica},
  volume={158},
  number={4},
  pages={750--763},
  year={2022},
  publisher={London Mathematical Society}
}

@article{arvidsson2023normality,
  title={Normality of minimal log canonical centers of threefolds in mixed and positive characteristic},
  author={Arvidsson, Emelie and Posva, Quentin},
  doi = {10.5802/aif.3780},
  volume={Online first},
  pages={33 p},
  journal={Annales de l'Institut Fourier},
  year={2026}
}

@article{arvidsson2023vanishing,
  title={On a vanishing theorem for birational morphisms of threefolds in positive and mixed characteristics},
  author={Arvidsson, Emelie},
  journal={arXiv preprint arXiv:2302.04420},
  year={2023}
}

@article{arvidsson2024properness,
  title={On the properness of the moduli space of stable surfaces over [1/30]},
  author={Arvidsson, Emelie and Bernasconi, Fabio and Patakfalvi, Zsolt},
  journal={Moduli},
  volume={1},
  pages={e3},
  year={2024},
  publisher={Cambridge University Press}
}

@article{bauer2002reduction,
  title={A reduction map for nef line bundles},
  author={Bauer, Thomas and Campana, Fr{\'e}d{\'e}ric and Eckl, Thomas and Kebekus, Stefan and Peternell, Thomas and Rams, S{\l}awomir and Szemberg, Tomasz and Wotzlaw, Lorenz},
  journal={Complex geometry},
  pages={27--36},
  year={2002},
  publisher={Springer}
}

@article{baum1979lefschetz,
  title={Lefschetz-riemann-roch for singular varieties},
  author={Paul Frank Baum and William Fulton and George Quart},
  journal={Acta Mathematica},
  year={1979},
  volume={143},
  pages={193-211},
  url={https://api.semanticscholar.org/CorpusID:122554018}
}

@article{bernasconi2024abundance,
  title={Abundance theorem for threefolds in mixed characteristic},
  author={Bernasconi, Fabio and Brivio, Iacopo and Stigant, Liam},
  journal={Mathematische Annalen},
  volume={388},
  number={1},
  pages={141--166},
  year={2024},
  publisher={Springer}
}

@article{bhatt2023globally,
  title={Globally-regular varieties and the minimal model program for threefolds in mixed characteristic},
  author={Bhatt, Bhargav and Ma, Linquan and Patakfalvi, Zsolt and Schwede, Karl and Tucker, Kevin and Waldron, Joe and Witaszek, Jakub},
  journal={Publications math{\'e}matiques de l'IH{\'E}S},
  volume={138},
  pages={1--159},
  year={2023},
  publisher={Springer}
}

@article{birkar2010existence,
  title={Existence of minimal models for varieties of log general type},
  author={Birkar, Caucher and Cascini, Paolo and Hacon, Christopher and McKernan, James},
  journal={Journal of the American Mathematical Society},
  volume={23},
  number={2},
  pages={405--468},
  year={2010}
}

@article{birkar2016existence,
  title={Existence of flips and minimal models for 3-folds in char $ p$},
  author={Birkar, Caucher},
  journal={Annales scientifiques de l'{\'E}cole normale sup{\'e}rieure},
  volume={49},
  number={1},
  pages={169--212},
  year={2016},
}

@article{birkar2017existence,
  title={Existence of Mori fibre spaces for 3-folds in char p},
  author={Birkar, Caucher and Waldron, Joe},
  journal={Advances in Mathematics},
  volume={313},
  pages={62--101},
  year={2017},
  publisher={Elsevier}
}

@article{cascini2015base,
 title = {On base point freeness in positive characteristic},
 author = {Cascini, Paolo and Tanaka, Hiromu and Xu, Chenyang},
 pages ={1239--1272},
 year = {2015},
 journal = {Annales Scientifiques de l'Ecole Normale Superieure},
 volume = {48},
 number = {5},
}

@article{das2019abundance,
  title={On the abundance problem for 3-folds in characteristic $ p> 5$},
  author={Das, Omprokash and Waldron, Joe},
  journal={Mathematische Zeitschrift},
  volume={292},
  number={3},
  pages={937--946},
  year={2019},
  publisher={Springer}
}

@article{das2022log,
author = {Das, Omprokash and Waldron, Joe},
title = {On the log minimal model program for threefolds over imperfect fields of characteristic $p>5$},
journal = {Journal of the London Mathematical Society},
volume = {106},
number = {4},
pages = {3895-3937},
doi = {https://doi.org/10.1112/jlms.12677},
url = {https://londmathsoc.onlinelibrary.wiley.com/doi/abs/10.1112/jlms.12677},
eprint = {https://londmathsoc.onlinelibrary.wiley.com/doi/pdf/10.1112/jlms.12677},
year = {2022}
}

@inproceedings{ekedahl1987foliations,
  title={Foliations and inseparable morphisms},
  author={Ekedahl, Torsten},
  booktitle={Algebraic geometry, Bowdoin, 1985},
  pages={139--149},
  year={1987}
}

@book{fulton2013intersection,
  title={Intersection theory},
  author={Fulton, William},
  volume={2},
  year={2013},
  publisher={Springer Science \& Business Media}
}

@article{gongyo2019rational,
  title={Rational points on log Fano threefolds over a finite field},
  author={Gongyo, Yoshinori and Nakamura, Yusuke and Tanaka, Hiromu},
  journal={Journal of the European Mathematical Society},
  volume={21},
  number={12},
  pages={3759--3795},
  year={2019}
}

@article{hacon2015three,
  title={On the three dimensional minimal model program in positive characteristic},
  author={Hacon, Christopher and Xu, Chenyang},
  journal={Journal of the American Mathematical Society},
  volume={28},
  number={3},
  pages={711--744},
  year={2015}
}

@article{hacon2022minimal,
author = {Christopher Hacon and Jakub Witaszek},
title = {{The minimal model program for threefolds in characteristic 5}},
volume = {171},
journal = {Duke Mathematical Journal},
number = {11},
publisher = {Duke University Press},
pages = {2193 -- 2231},
year = {2022},
doi = {10.1215/00127094-2022-0024},
URL = {https://doi.org/10.1215/00127094-2022-0024}
}

@article{hacon2022relative,
  title={On the relative minimal model program for threefolds in low characteristics},
  author={Hacon, Christopher and Witaszek, Jakub},
  journal={Peking Mathematical Journal},
  volume={5},
  number={2},
  pages={365--382},
  year={2022},
  publisher={Springer}
}

@article{hashizume2020minimal,
  title={Minimal model program for log canonical threefolds in positive characteristic},
  author={Kenta Hashizume and Yusuke Nakamura and Hiromu Tanaka},
  journal={Mathematical Research Letters},
  volume={27},
  pages={1003-–1054},
  year={2017},
}

@article{kawamata1992abundance,
  title={Abundance theorem for minimal threefolds},
  author={Kawamata, Yujiro},
  journal={Inventiones mathematicae},
  volume={108},
  pages={229--246},
  year={1992},
  publisher={Springer}
}

@article{keel1994log,
  title={Log abundance theorem for threefolds},
  author={Keel, Sean and Matsuki, Kenji and McKernan, James},
  journal={Duke Mathematical Journal},
  volume={75},
  number={1},
  pages={99--119},
  year={1994},
  publisher={Duke University Press}
}

@article{keeler2003ample,
  title={Ample filters of invertible sheaves},
  author={Keeler, Dennis S},
  journal={Journal of Algebra},
  volume={259},
  number={1},
  pages={243--283},
  year={2003},
  publisher={Elsevier}
}

@article{kleiman1969geometry,
  title={Geometry on Grassmannians and applications to splitting bundles and smoothing cycles},
  author={Kleiman, Steven L},
  journal={Publications Math{\'e}matiques de l'IH{\'E}S},
  volume={36},
  pages={281--297},
  year={1969}
}

@inproceedings{kollar1992flips,
  title={Flips and abundance for algebraic threefolds},
  author={Koll{\'a}r, J{\'a}nos},
  booktitle={Lectures in Utah Summer Seminar 1991, Asterisque},
  volume={211},
  year={1992}
}

@book{kollar1998birational,
  title={Birational geometry of algebraic varieties},
  author={Koll{\'a}r, Janos and Mori, Shigefumi and Clemens, Charles Herbert and Corti, Alessio},
  volume={134},
  year={1998},
  publisher={Cambridge university press}
}

@book{kollar2013singularities,
  title={Singularities of the minimal model program},
  author={Koll{\'a}r, J{\'a}nos},
  volume={200},
  year={2013},
  publisher={Cambridge University Press}
}

@article{langer2000chern,
  title={Chern classes of reflexive sheaves on normal surfaces},
  author={Langer, Adrian},
  journal={Mathematische Zeitschrift},
  volume={235},
  number={3},
  pages={591--614},
  year={2000},
  publisher={Springer}
}

@article{langer2004semistable,
  title={Semistable sheaves in positive characteristic},
  author={Langer, Adrian},
  journal={Annals of mathematics},
  volume={159},
  pages={251--276},
  year={2004}
}

@article{langer2015generic,
  title={Generic positivity and foliations in positive characteristic},
  author={Langer, Adrian},
  journal={Advances in Mathematics},
  volume={277},
  pages={1--23},
  year={2015},
  publisher={Elsevier}
}

@article{langer2024bogomolov,
  title={Bogomolov’s inequality and Higgs sheaves on normal varieties in positive characteristic},
  author={Langer, Adrian},
  journal={Journal f{\"u}r die reine und angewandte Mathematik (Crelles Journal)},
  volume={2024},
  number={810},
  pages={1--48},
  year={2024},
  publisher={De Gruyter}
}

@article{langer2025intersection,
  title={Intersection theory and Chern classes on normal varieties},
  author={Langer, Adrian},
  journal={Journal of the London Mathematical Society},
  volume={112},
  number={1},
  pages={e70244},
  year={2025},
  publisher={Wiley Online Library}
}

@book{lazarsfeld2017positivity,
  title={Positivity in algebraic geometry I: Classical setting: line bundles and linear series},
  author={Lazarsfeld, Robert K},
  volume={48},
  year={2017},
  publisher={Springer}
}

@article{liu2016geography,
  title={Geography of Gorenstein stable log surfaces},
  author={Liu, Wenfei and Rollenske, S{\"o}nke},
  journal={Transactions of the American Mathematical Society},
  volume={368},
  number={4},
  pages={2563--2588},
  year={2016},
  publisher={American Mathematical Society}
}

@article{luisa1998axiomatic,
  title={Axiomatic theory for transversality and Bertini type theorems},
  author={Luisa Spreafico, Maria},
  journal={Archiv der Mathematik},
  volume={70},
  pages={407--424},
  year={1998},
  publisher={Springer}
}

@inproceedings{mehta2006homogeneous,
  title={Homogeneous bundles in characteristic p},
  author={Mehta, VB and Ramanathan, A},
  booktitle={Algebraic Geometry—Open Problems: Proceedings of the Conference Held in Ravello, May 31--June 5, 1982},
  pages={315--320},
  year={1983},
  volumn={997},
  organization={Lecture Notes in Mathematics, vol 997, Springer Berlin Heidelberg}
}

@article{posva2023abundance,
  title={Abundance for slc surfaces over arbitrary fields},
  author={Posva, Quentin},
  journal={{\'E}pijournal de G{\'e}om{\'e}trie Alg{\'e}brique},
  volume={7},
  pages={1--22},
  year={2023},
  publisher={Episciences. org}
}

@article{reid1985young,
  title={Young person’s guide to canonical singularities},
  author={Reid, Miles},
  journal={Algebraic geometry, Bowdoin},
  volume={46},
  pages={345--414},
  year={1985}
}

@article{sato2020general,
  title={General hyperplane sections of threefolds in positive characteristic},
  author={Sato, Kenta and Takagi, Shunsuke},
  journal={Journal of the Institute of Mathematics of Jussieu},
  volume={19},
  number={2},
  pages={647--661},
  year={2020},
  publisher={Cambridge University Press}
}

@article{sato2025general,
  title={General hyperplane sections of log canonical threefolds in positive characteristic},
  author={Sato, Kenta},
  journal={Journal of the Institute of Mathematics of Jussieu},
  volume={24},
  number={5},
  pages={1867--1894},
  year={2025},
  publisher={Cambridge University Press}
}

@book{serre1977linear,
  title={Linear representations of finite groups},
  author={Serre, Jean-Pierre},
  volume={42},
  year={1977},
  publisher={Springer}
}

@misc{stacks-project,
  author       = {The {Stacks project authors}},
  title        = {The Stacks project},
  howpublished = {\url{https://stacks.math.columbia.edu}},
  year         = {2024},
}

@article{takamatsu2023m,
  title={Minimal model program for semi-stable threefolds in mixed characteristic},
  author={Takamatsu, Teppei and Yoshikawa, Shou},
  journal={Journal of Algebraic Geometry},
  volume={32},
  pages={429--476},
  year={2023}
}

@article{tanaka2012x,
  title={The {X}-method for klt surfaces in positive characteristic},
  author={Tanaka, Hiromu},
  journal={Journal of Algebraic Geometry},
  volume={24},
  pages={605--628},
  year={2015}
}

@article{waldron2017finite,
title = "Finite generation of the log canonical ring for 3-folds in char p",
author = "Joe Waldron",
year = "2017",
doi = "10.4310/MRL.2017.v24.n3.a14",
language = "English (US)",
volume = "24",
pages = "933--946",
journal = "Mathematical Research Letters",
issn = "1073-2780",
publisher = "International Press of Boston, Inc.",
number = "3",
}

@article{waldron2018lmmp,
  title={The {LMMP} for log canonical 3-folds in characteristic $p>5$},
  author={Waldron, Joe},
  journal={Nagoya Mathematical Journal},
  volume={230},
  pages={48--71},
  year={2018},
  publisher={Cambridge University Press}
}

@article{waldron2023mori,
  title={Mori fibre spaces for $3 $-folds over imperfect fields},
  author={Waldron, Joe},
  journal={arXiv preprint arXiv:2303.00615},
  year={2023}
}

@article{witaszek2021canonical,
  title={On the canonical bundle formula and log abundance in positive characteristic},
  author={Witaszek, Jakub},
  journal={Mathematische Annalen},
  volume={381},
  number={3},
  pages={1309--1344},
  year={2021},
  publisher={Springer}
}

@article{xu2019nonvanishing,
author = {Chenyang Xu and Lei Zhang},
title = {{Nonvanishing for $3$-folds in characteristic $p> 5$}},
volume = {168},
journal = {Duke Mathematical Journal},
number = {7},
publisher = {Duke University Press},
pages = {1269 -- 1301},
year = {2019},
doi = {10.1215/00127094-2018-0066},
URL = {https://doi.org/10.1215/00127094-2018-0066}
}

@article{xu2024note,
  title={NOTE ON THE THREE-DIMENSIONAL LOG CANONICAL ABUNDANCE IN CHARACTERISTIC $>3$},
  author={Xu, Zheng},
  journal={Nagoya Mathematical Journal},
  doi={10.1017/nmj.2024.3},
  volume={255},
  pages={694--723},
  year={2024}
}

@article{zhang2019abundance,
  title={Abundance for non-uniruled 3-folds with non-trivial Albanese maps in positive characteristics},
  author={Zhang, Lei},
  journal={Journal of the London Mathematical Society},
  volume={99},
  number={2},
  pages={332--348},
  year={2019},
  publisher={Wiley Online Library}
}

@article{zhang2020abundance,
  title={Abundance for 3-folds with non-trivial Albanese maps in positive characteristic},
  author={Zhang, Lei},
  journal={Journal of the European Mathematical Society},
  volume={22},
  number={9},
  pages={2777--2820},
  year={2020}
}

@article{zhang2023frobenius,
  title={Frobenius stable pluricanonical systems on threefolds of general type in positive characteristic},
  author={Zhang, Lei},
  journal={Algebra \& Number Theory},
  volume={16},
  number={10},
  pages={2339--2384},
  year={2023},
  publisher={Mathematical Sciences Publishers}
}

\end{document}